\documentclass[%
 aip,
 amsmath,amssymb,
reprint, onecolumn %%%%% Comment onecolumn to override one column format %%%%%%%
]{revtex4-1}

\usepackage{graphicx}% Include figure files
\usepackage{subcaption}% Include subfigures
\usepackage{dcolumn}% Align table columns on decimal point
\usepackage{bm}% bold math
\usepackage[utf8]{inputenc}
\usepackage[T1]{fontenc}
\usepackage{mathptmx}
\usepackage{etoolbox}
\usepackage{comment}
\usepackage{color}
\usepackage{fancyhdr}
\usepackage{footmisc}
\usepackage{hyperref}
\usepackage[margin=22mm]{geometry}
\usepackage[UKenglish]{datetime}
\usepackage{graphicx}
\usepackage{wrapfig}
\usepackage{enumerate}
\usepackage{tikz-cd}

\newcommand{\gaw}[1]{\textcolor{purple}{#1}}

\graphicspath{{./img/}}

\makeatletter
\def\@email#1#2{%
 \endgroup
 \patchcmd{\titleblock@produce}
  {\frontmatter@RRAPformat}
  {\frontmatter@RRAPformat{\produce@RRAP{*#1\href{mailto:#2}{#2}}}\frontmatter@RRAPformat}
  {}{}
}%
\makeatother
\begin{document}

\preprint{AIP/123-QED}

\title[High order interpolation of magnetic fields with vector potential reconstruction for particle simulations]{High order interpolation of magnetic fields with vector potential reconstruction for particle simulations}
% Force line breaks with \\
\author{O. Beznosov}
\affiliation{%
Los Alamos National Laboratory
}%Lines break
  \email{obeznosov@lanl.gov.}
  %automatically or can be forced with \\
 \author{J. Bonilla}%
 % \email{Second.Author@institution.edu.}
\affiliation{%
Los Alamos National Laboratory
}%Lines break
% \affiliation{
% Authors' institution and/or address%\\This line break forced with \textbackslash\textbackslash
% }%
\author{G. A. Wimmer}%
 % \email{Second.Author@institution.edu.}
\affiliation{%
Los Alamos National Laboratory
}%%Lines break
\author{X.-Z. Tang}
\affiliation{%
Los Alamos National Laboratory
}%

\date{\today}% It is always \today, today,

             \begin{abstract}
              We propose a method for interpolating divergence-free continuous magnetic fields via vector potential reconstruction using Hermite interpolation, which ensures high-order continuity for applications requiring adaptive, high-order ordinary differential equation (ODE) integrators, such as the Dormand-Prince method. The method provides \(C(m)\) continuity and achieves high-order accuracy, making it particularly suited for particle trajectory integration and Poincaré section analysis under optimal integration order and timestep adjustments. Through numerical experiments, we demonstrate that the Hermite interpolation method preserves volume and continuity, which are critical for conserving toroidal canonical momentum and magnetic moment in guiding center simulations, especially over long-term trajectory integration. Furthermore, we analyze the impact of insufficient derivative continuity on Runge-Kutta schemes and show how it degrades accuracy at low error tolerances, introducing discontinuity-induced truncation errors. Finally, we demonstrate performant Poincaré section analysis in two relevant settings of field data collocated from finite element meshes.
             \end{abstract}

\maketitle

\section{\label{sec:Introduction}Introduction}

Electromagnetics and plasma simulations are subject to a fundamental physics constraint on the magnetic field \(\mathbf{B}\), which is zero divergence (\(\nabla\cdot\mathbf{B}=0\)) or the absence of magnetic monopoles. There are two common approaches to address this physics constraint in discrete mathematics. One less deliberate approach is to control the non-vanishing \(\nabla\cdot\mathbf{B}\) below the level of discretization error, by introducing a so-called divergence cleaning term in the equation.~\cite{Zhang-Feng-Frontiers-2016,TOTH2000605,DEDNER2002645}

More elaborate approaches aim to enforce \(\nabla\cdot\mathbf{B}=0\) in the design of numerical discretizations. A representative method is the mimetic finite difference method,~\cite{LIPNIKOV20141163,LipnikovMoultonSvyatskiy2008,BeiraodaVeigaGyryaLipnikovManzini2009,BeiraodaVeigaGyryaLipnikovManzini2009,BeiraodaVeiga2010,BeiraodaVeigaDroniouManzini2011,CangianiManziniRusso2009,BrezziLipnikovShashkovSimoncini2007} in which the magnetic field is discretized as a surface normal on each cell surface, shared by neighboring cells on both sides of the cell surface. The discrete numerical divergence, evaluated at the cell center, strictly remains divergence-free if the initial field satisfies \(\nabla\cdot\mathbf{B}=0\). The Yee mesh in particle-in-cell simulation~\cite{Yee-IEEE-1966} is an early realization of this approach and has found applications in a wide range of plasma physics problems.

\begin{comment}
Electromagnetics and plasma simulations have a fundamental physics
constraint on the magnetic field $\mathbf{B},$ which is zero
divergence ($\nabla\cdot\mathbf{B}=0$) or absence of magnetic
monopoles.  There are two common approaches to deal with this physics
constraint in discrete mathematics.  A less deliberate approach is to
control the non-vanishing $\nabla\cdot\mathbf{B}$ below the level of
discretization error, by introducing a so-called divergence cleaning
term in the
equation.~\cite{Zhang-Feng-Frontiers-2016,TOTH2000605,DEDNER2002645}
More elaborate approaches aim to impose $\nabla\cdot\mathbf{B}=0$ in
the design of numerical discretization.  A representative method is
mimetic finite
difference,\cite{LIPNIKOV20141163,LipnikovMoultonSvyatskiy2008,BeiraodaVeigaGyryaLipnikovManzini2009,BeiraodaVeigaGyryaLipnikovManzini2009,BeiraodaVeiga2010,BeiraodaVeigaDroniouManzini2011,CangianiManziniRusso2009,BrezziLipnikovShashkovSimoncini2007}
in which the magnetic field is discretized as a surface normal on each
cell surface, shared by neighboring cells on both sides of the cell
surface. The discrete numerical divergence, which is evaluated at the
cell center, remains strictly divergence-free if the initial field
satisfies $\nabla\cdot\mathbf{B}=0.$ The Yee mesh in particle-in-cell
simulation~\cite{Yee-IEEE-1966} is an early realization of this
approach and has found applications in a wide-range of plasma physics
problems.
\end{comment}
For integrating magnetic field line characteristics, which is required in diagnosing magnetic topology and/or evolving parallel transport, and for tracing charged particle characteristics, especially that of the guiding center, one would desire \(\nabla\cdot\mathbf{B}=0\) everywhere in the simulation domain. Standard interpolation schemes of a discrete mimetic \(\mathbf{B}\) field cannot provide this property, even though they are commonly used in Poincaré section analysis and particle orbit tracing. The situation is further aggravated if \(\nabla\cdot\mathbf{B}\) is not even satisfied discretely in the numerical discretization, in which case one often finds Poincaré section analysis impossible to perform or that the results are simply nonsensical. Finn and Chacon~\cite{finn2005} proposed an approach that enforces \(\nabla\cdot\mathbf{B}=0\) through reconstruction. The idea is to first construct the vector potential \(\mathbf{A}_D\) on discrete grid points from the numerical/discrete \(\mathbf{B}_D\), and then spline fit the discrete \(\mathbf{A}_D\) for a globally smooth vector potential \(\mathbf{A}\). The final step involves taking the analytical derivative of the vector potential components for \(\mathbf{B} = \nabla\times\mathbf{A}\), which is divergence-free everywhere in the simulation domain. This approach has proven to produce sensible Poincaré section plots; namely, nested flux surfaces remain sharp, and separatrix crossing can be clearly identified in the transition to field line chaos. In other words, the reconstruction scheme can reliably project a discrete vector field, even if it is not divergence-free, onto the space of divergence-free vectors.

The purpose of this paper is to quantify the accuracy of such a reconstruction scheme. This appears to be a meaningful problem only if the discrete \(\mathbf{B}_D\) field is mimetic in the first place. In addition to Poincaré section analysis, we are particularly interested in how such a reconstructed magnetic field impacts the accuracy of the guiding center integration. This is motivated by the application of Lagrangian and semi-Lagrangian integration of the relativistic drift-kinetic equation, which has found critical importance in understanding the runaway electron dynamics in a tokamak disruption.~\cite{mcdevitt2019avalanche,mcdevitt_pre_2023} 

There are a number of factors that determine the accuracy of the reconstruction scheme. The first concerns the numerical integration that constructs the vector potential from the magnetic field data on a computational grid. The second is the interpolation scheme for the discrete vector potential data that produces the globally smooth \(\mathbf{A}\), from which \(\mathbf{B} = \nabla\times\mathbf{A}\) can be evaluated everywhere in the simulation domain. The third is the numerical integration of the guiding center trajectory using the interpolated field data. 

In place of the cubic spline, we are interested in how Hermite interpolation can provide a useful alternative. This was partly motivated by the recent success in applying Hermite interpolation to construct highly accurate and stable solvers of partial differential equations (PDEs)~\cite{goodrich2006hermite, ApploHermite18, beznosov2021hermite}, and by the relative ease in computing the high-order Hermite interpolants, compared with the cubic spline. We were also attracted to Hermite interpolation for its attribute of not overfitting. The basics of Hermite interpolation are explained in section~\ref{ssec:Hermite_interpolation}.

What we have learned is that achieving accuracy and smoothness of the interpolated magnetic field requires care in the very first step of reconstructing the vector potential \(\mathbf{A}\) from the discrete \(\mathbf{B}_D\) data on the grid. This concerns the spatial integration of the magnetic field for the vector potential components, as discussed in section~\ref{sec:Flux_reconstruction}. For accuracy, we find that analytical integration of Hermite interpolation of the \(\mathbf{B}\) data is a superior approach compared with the usual line integration by numerical quadratures. Since Hermite interpolation, by matching the derivatives up to the \(m\)th order with the data, can achieve an order of \(2m + 1\) in the interpolating polynomials, we anticipate globally smooth fitting of the vector potential, as long as the numerical derivatives of the discrete field data on the grid points are reasonably well-behaved. 

The smoothness in the reconstructed \(\mathbf{A}\) and the resulting divergence-free magnetic field \(\mathbf{B}\) is conducive to higher-order integration of the guiding center trajectories. This is because high-order ordinary differential equation (ODE) integration schemes, such as Runge-Kutta (RK), are based on Taylor expansion. Thus, their accuracy depends on the existence of smooth high-order derivatives. Since the smoothness of the trajectories follows from the smoothness of the governing field~\cite{lekien2005}, for coherent and consistent reconstruction of particle trajectories, the high-order smoothness of the background fields becomes paramount.

\begin{comment}
What we have learned is that achieving accuracy and smoothness of the
interpolated magnetic field requires care in the very first step of
reconstructing the vector potential $\mathbf{A}$ from the discrete
$\mathbf{B}_D$ data on the grid. This concerns the spatial integration
of the magnetic field for the vector potential components,
section~\ref{sec:Flux_reconstruction}. For accuracy, we find that
analytical integration of Hermite interpolation of the $\mathbf{B}$
data is a superior approach, compared with the usual line integration
by numerical quadratures.  Since by matching the derivatives up to the
$m$th order with the data, Hermite interpolation can achieve the order
of $2m + 1$ in the interpolating polynomials, and hence we anticipate globally
smooth fitting of the vector potential, as long as the numerical
derivatives of the discrete field data on the grid points are
reasonably behaved.  The smoothness in the reconstructed $\mathbf{A}$ and
the resulting divergence-free magnetic field $\mathbf{B},$ is
conducive to higher order integration of the guiding center
trajectories, since high order ordinary differential equation (ODE) integration scheme, such as
Runge-Kutta (RK), is based on Taylor expansion, so its accuracy depends on
the existence of smooth high-order derivatives. Since the
smoothness of the trajectories follows from the smoothness of the
governing field \cite{lekien2005}, for coherent and
consistent reconstruction of particle trajectories the high order
smoothness of the background fields becomes paramount.
\end{comment}
Our numerical experiments indicate that, indeed, whether the RK method can optimally converge for Poincaré section analysis and particle trajectory integration depends on the smoothness of the reconstructed magnetic field. As theoretically shown~\cite{hairer1993}, the order conditions must be satisfied for the RK method, say of the 5th order as in the Dormand-Prince 4(5) pair, to properly adjust the timestep. This requires derivatives up to the fourth order to cancel in certain linear combinations of RK stages, defined by a Butcher tableau.

The rest of the paper is organized as follows. In section~\ref{sec:OrderReduction}, we demonstrate the reduction in the order of accuracy of Runge-Kutta methods on a simple problem to further motivate the use of Hermite interpolation and show what we expect to see in our numerical experiments with magnetic fields. Section~\ref{ssec:Hermite_interpolation} introduces the Hermite interpolation method. The construction of a vector potential from the discrete magnetic field data is presented in section~\ref{sec:Flux_reconstruction}. Various numerical experiments, including the use of magnetic field data from both finite element and mimetic finite difference MHD codes, are reported in section~\ref{sec:Numerical_experiments}. Section~\ref{sec:Summary} concludes the paper with a summary.

\section{\label{sec:OrderReduction}Convergence order reduction in high order ODE solvers}

As conventionally defined\cite{hairer1993}, an RK method with
time step $\Delta t$ has order $p$ if, for sufficiently smooth problems,
the Taylor series for the exact solution $y(t_0 + \Delta t)$ and
numerical solution $y_1$ agree with each other up to the term
$\Delta t^p$. In practice, Runge-Kutta type methods are also 
frequently used to solve ODEs whose right-hand sides do not provide
such smoothness or in some cases are even discontinuous. This may
for instance occur for ODEs resulting from PDEs that are discretized
in space using non-smooth finite element or finite difference methods.
Correspondingly, the lack of a Taylor series may then deteriorate the
expected convergence rate, which we will demonstrate with a simple
example of a non-autonomous ODE. However, smoothness can be
restored by means of Hermite interpolation, which is particularly useful
where high-order ODE integration is desirable.

\begin{figure}[ht]
\begin{center}
\includegraphics[width=.3\textwidth]{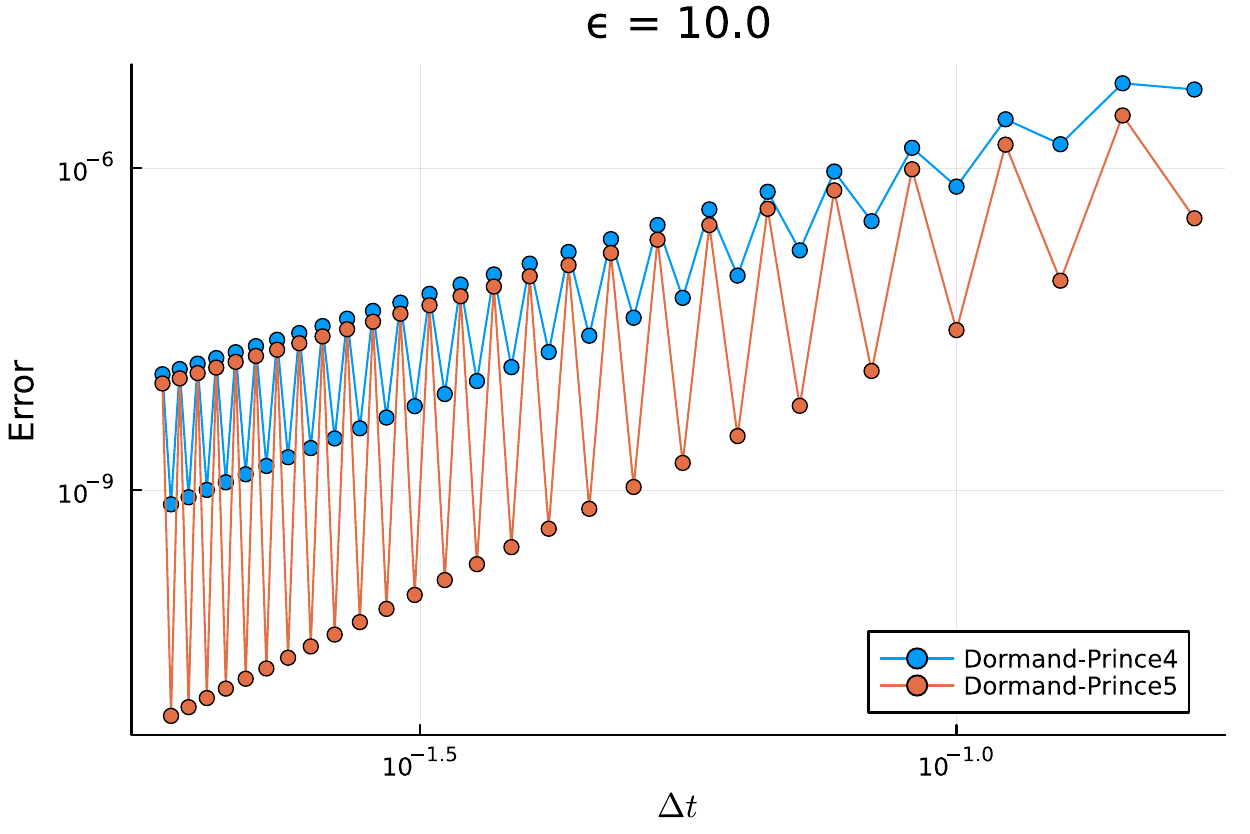}
\includegraphics[width=.3\textwidth]{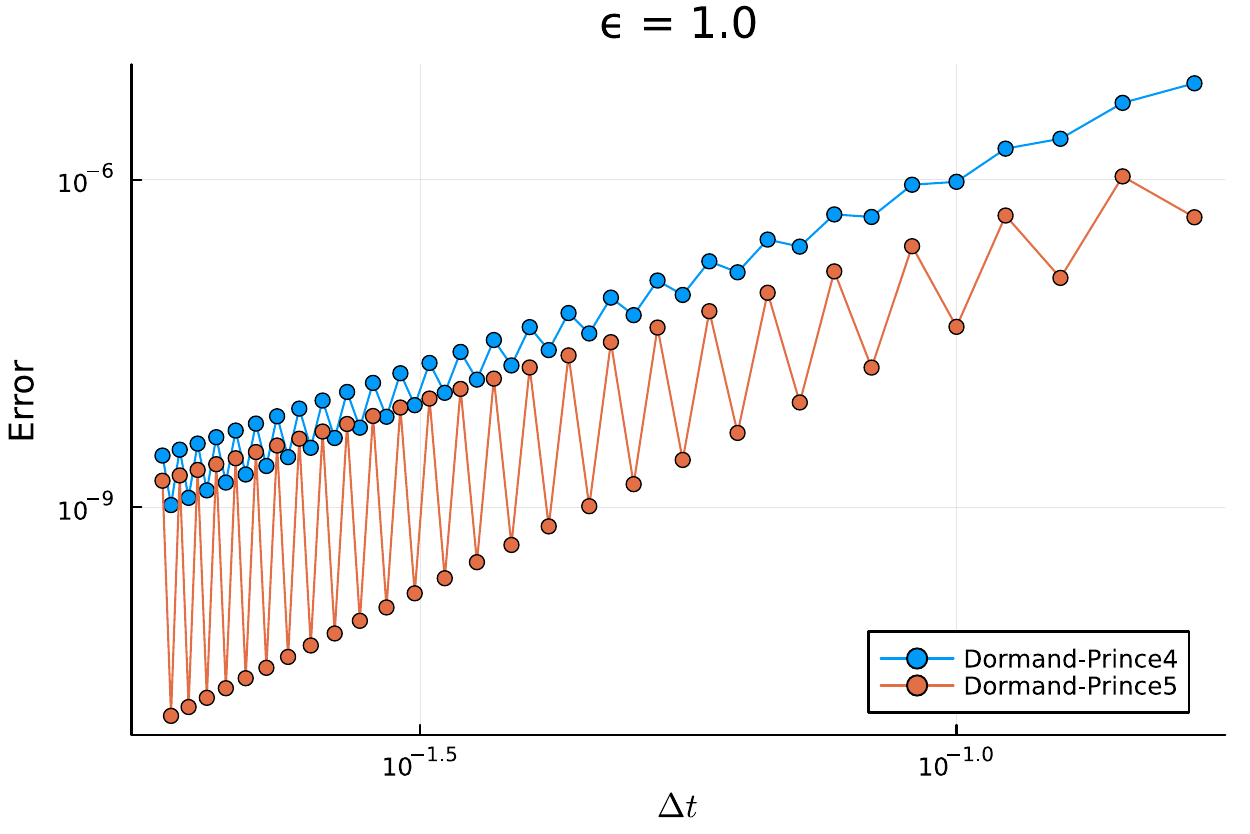}
\includegraphics[width=.3\textwidth]{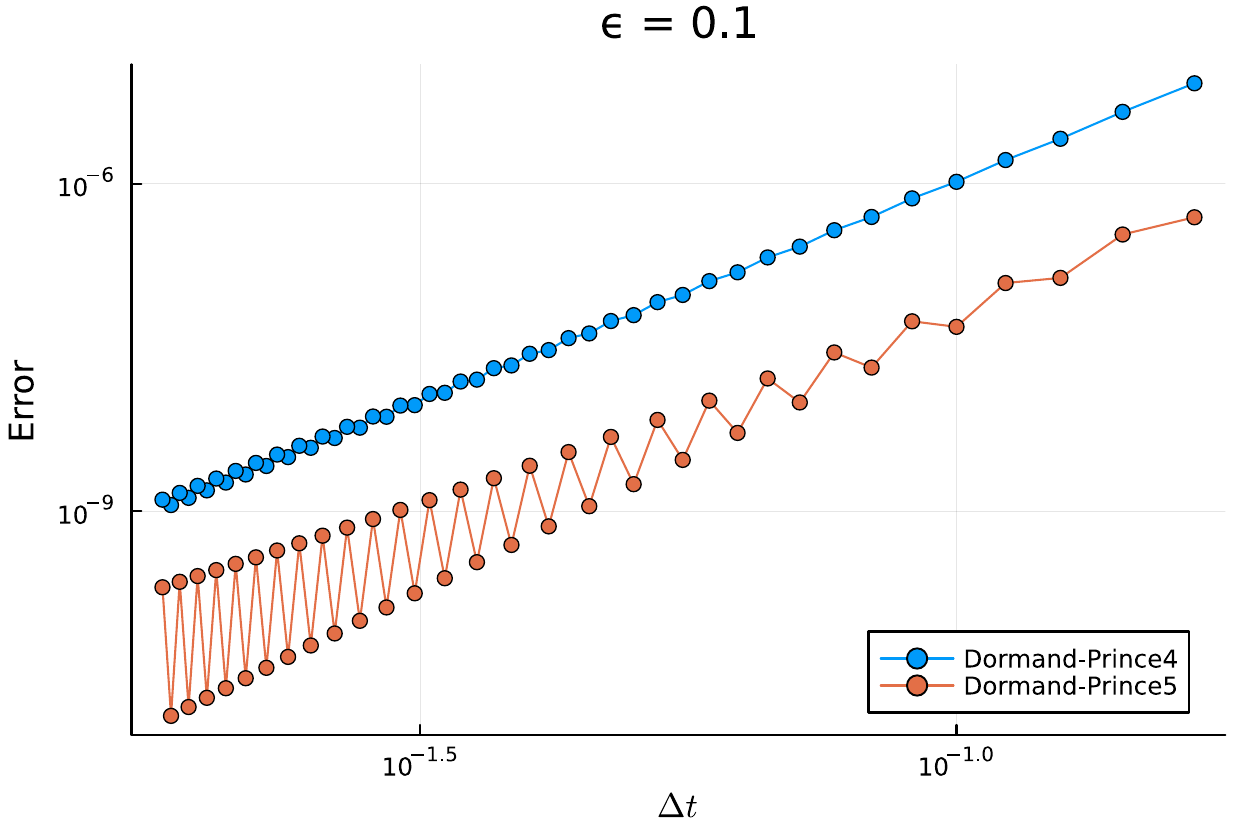}
\includegraphics[width=.3\textwidth]{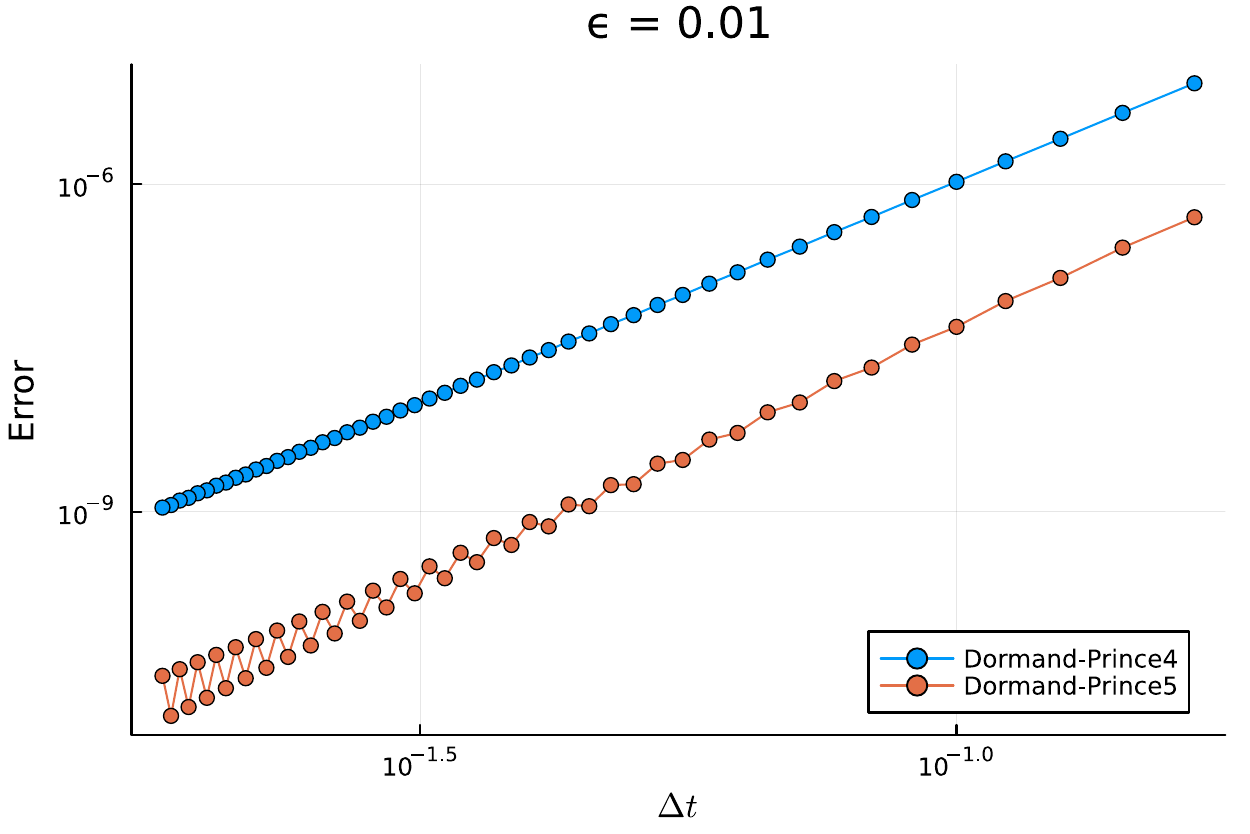}
\includegraphics[width=.3\textwidth]{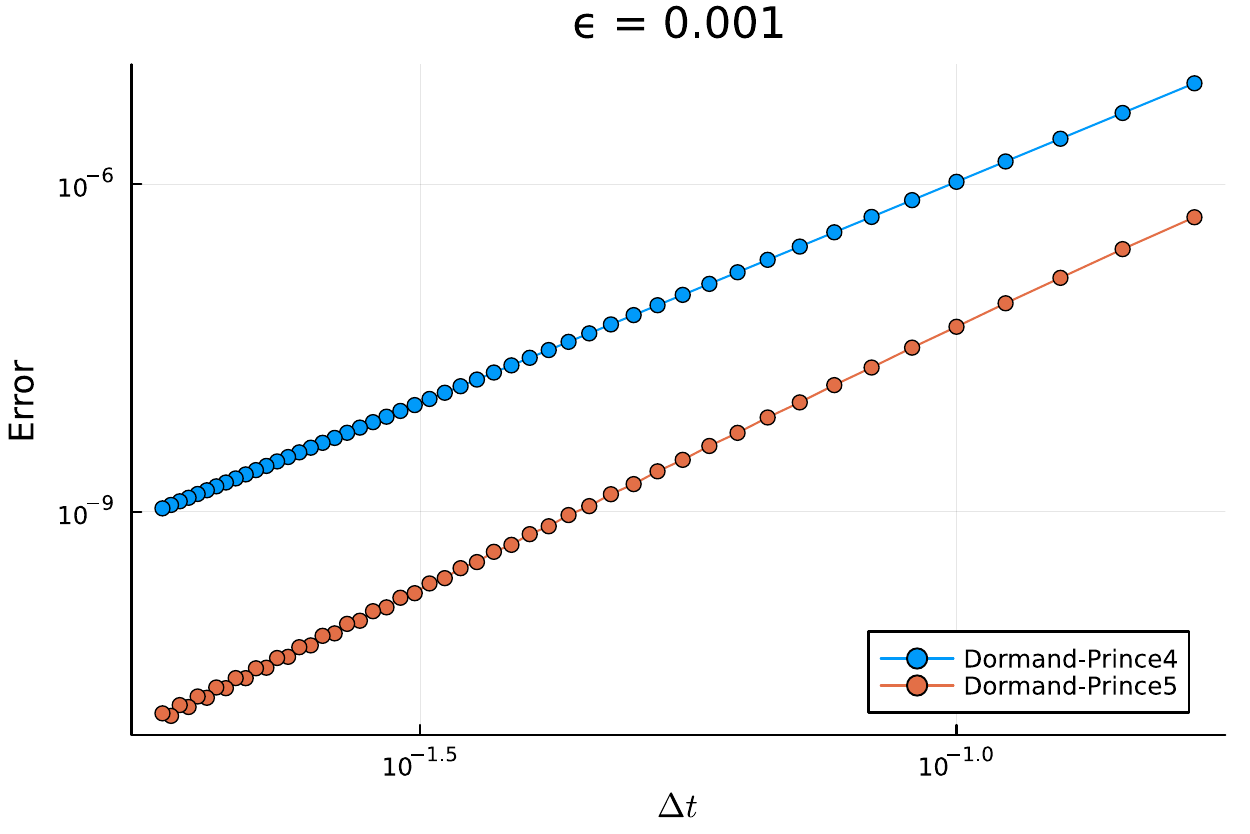}
\includegraphics[width=.3\textwidth]{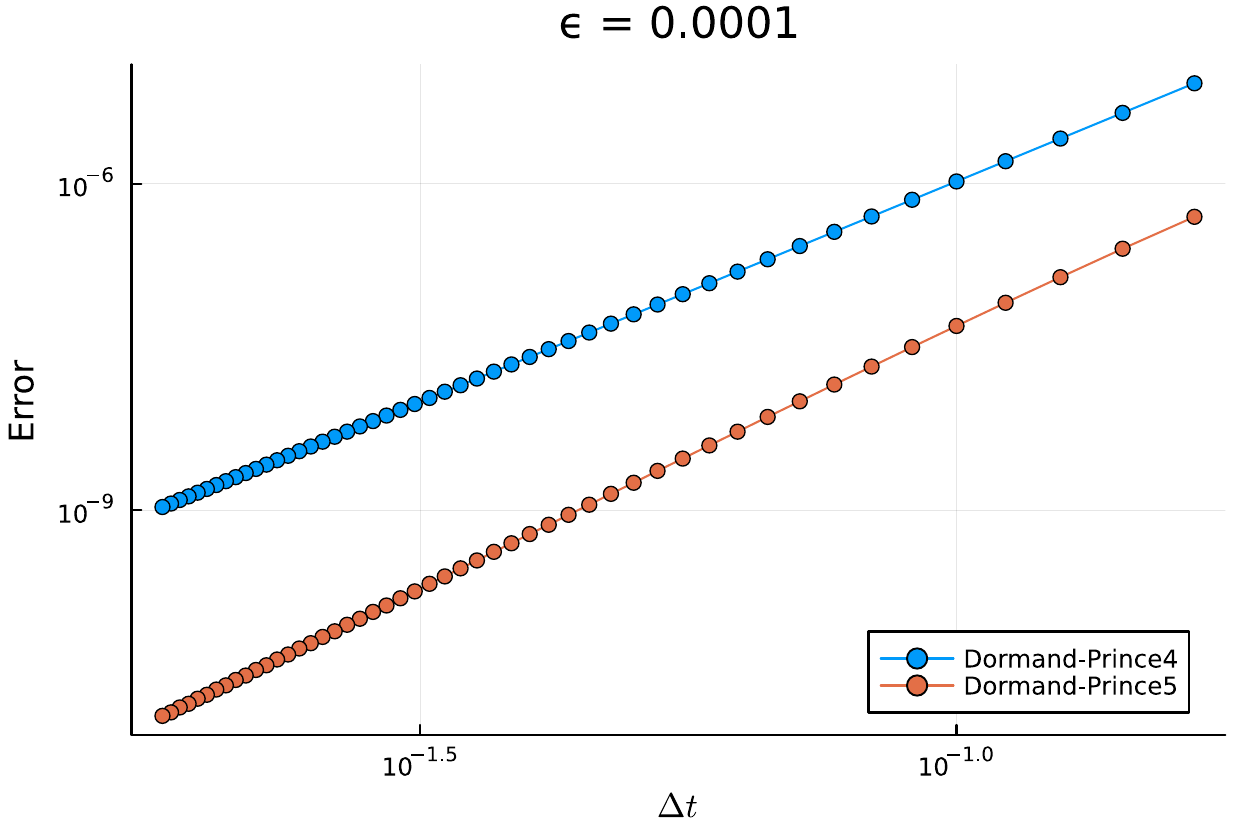}
\end{center}
\caption{\label{fig:oderk} Convergence rates for 4th order (blue curves) and 5th order (red curves) solutions to the Dormand-Prince Runge-Kutta method. Curves display numerical error from integrating \eqref{eq:ode} for a right-hand side with second derivative discontinuity of varying magnitude $\epsilon \in [10^{-4}, 10]$. The error is shown as a function of time step size, and oscillations arise from passing through time step sizes that alternatingly lead an even and odd number of total time steps.}
\end{figure}

Consider a first order ODE,
\begin{equation}
  y' = f(t),\ \ y(0) =0, \label{eq:ode}
\end{equation}
with a piecewise right-hand side
\begin{equation}
  f(t) = \begin{cases}
    f_1(t), \ \ &{t \le \frac 12}, \\
    f_2(t), \ \ &{t > \frac 12}.
  \end{cases}
\end{equation}
The exact solution at the end point $t = 1$ can be written as an integral
over the full time interval,
\begin{equation}
y(1) = \int_0^1 f(t)\, {\rm d}t.
\end{equation}
To further study the aforementioned effect of non-smoothness, we set
up a numerical experiment choosing $f_1$ and $f_2$ such that
$f_1(\frac 12) = f_2(\frac 12)$, $f '_1(\frac 12) = f'_2(\frac 12)$
and ${f''_1}(\frac 12) - {f''_2}(\frac 12) = \epsilon$. Numerically integrating
\eqref{eq:ode} with a fixed time step using the embedded Dormand-Prince
Runge-Kutta pair\cite{dormand1980family}, we obtain convergence plots as
given in FIG.~\ref{fig:oderk}, depicting both the method's $5^{th}$ and
embedded $4^{th}$ order solutions for a range of values for $\epsilon$.
For an even total number of time steps, both solutions converge at the
expected rate, while for an odd number the solutions' convergence
order deteriorates from $5$ to $3$ and from $4$ to $3.5$ for the two
solutions, respectively. The decreased convergence rate for an odd
number of time steps is in line with our above discussion on smoothness.
On the other hand, the optimal convergence for an even number of time 
steps follows since in this case one of the discrete points in time is exactly at the 2nd
derivative discontinuity at $t = \frac 12$, and therefore locally, the right-hand
side function fed into the RK algorithm is always perfectly smooth. Finally,
the negative effect for an odd number of time steps decreases as $\epsilon$
is decreased, that is as the discontinous gap at $t= \frac{1}{2}$ narrows. In
particular, small values of $\epsilon$ -- corresponding to a ``almost regular''
right-hand side --  lead to convergence rates that may still be deemed
acceptable in practice.

Similar cases can be constructed for autonomous systems, but for the
general case $y'=f(y)$, standard time stepping methods are unlikely
to always step exactly on given discontinuities. Importantly, as seen in
the above example for the Dormand-Prince method, adaptive
timestepping capabilities may deteriorate due to the reduced convergence
rates in the presence of non-smoothness, forcing the method to march at a
minimal timestep\cite{hairer1993}. To avoid this, one could employ 
dense output -- using e.g. Hermite interpolation -- and trajectory inversion
techniques in order to not step over discontinuities. Alternatively, one could
impose high-order continuity to begin with by applying a high-order Hermite
interpolation method on the ODE's right-hand side, and this is the approach
covered in this work.

\section{\label{ssec:Hermite_interpolation}Hermite interpolation}

We briefly present the Hermite interpolation method here and refer the reader
to the original work on Hermite methods for wave equations\cite{ApploHermite18} for error analysis and other details. For ease of exposition,
we first demonstrate a one-dimensional approximation and then show how to
extend it to higher dimensions.

Let the spatial domain be $D = [a,b]$, discretized by a primal grid
\begin{equation}
  \Omega_p = \left\{
    x_i:\ \ x_i = a + ih,\ \ i = 0,\dots,N
    \right\},
\end{equation}
with corresponding dual grid
\begin{equation}
  \Omega_d = \left\{x_j:\ \ x_j = a + jh,\ \ j = \frac 12,\dots, N - \frac 12\right\}.
\end{equation}

At the primal grid $\Omega_p$, we store the target function $f$ and its $m$
derivatives, represented by Taylor polynomials of order $m$. At the dual
grid $\Omega_d$, we center local polynomials that match the target
function and its $m$ derivatives at primal grid points.
\begin{comment}
  At each grid point, the approximation to the solution is represented by its degrees of freedom (DOF) that approximate the values and spatial derivatives of target function $f$. Equivalently, the approximation $f$ can be represented as polynomials centered at grid points. The Taylor coefficients of these polynomials are scaled versions of the DOF. \textcolor{blue} {To achieve the order $(2m+1)$ of the interpolating polynomials and continuity of $m$ first derivatives we require the }\textcolor{blue}{$m$ first derivatives} of $f$ to be stored at each grid point.
\end{comment}
First, we obtain the polynomials as approximations to the target function on the primal grid
\begin{equation}
f(x) \approx \sum_{l = 0}^{m} \hat u_{l} \left(\frac{x-x_i}h \right)^l,\ \  x_i \in \Omega_p. \label{coef1D}
\end{equation}

The coefficients $\hat u_l$ are assumed to be accurate approximations to the scaled Taylor coefficients of the target function. If expressions for the derivatives were known, for example from an underlying high order method, we could simply set
\begin{equation}
  \hat u_{l} = \frac{h^l}{l!}\frac{{\rm d}^l f}{{\rm d}x^l}(x_i), \ \ x_i \in \Omega_p.
\end{equation}
If the derivatives are not available, we can instead use finite difference approximations on a finer grid.
In our numerical experiments, these coefficients are evaluated with centered finite differences\cite{fornberg1988generation} as follows
\begin{align}
  \hat u_{0} &= f(x_i), \label{eq:fd1}\\
  \hat u_{1} &= h \cdot \frac{-f(x_i + 2\Delta x) + 8f(x_i+ \Delta x) - 8f(x_i - \Delta x) + f(x_i - 2\Delta x)} {12 \Delta x}, \label{eq:fd2}\\
  \hat u_{2} &= \frac{h^2}2 \cdot \frac{-f(x_i + 2\Delta x)  + 16f(x_i+ \Delta x) - 30f(x_i) + 16f(x_i - \Delta x) - f(x_i - 2\Delta x)}{12 \Delta x^2 }, \label{eq:fd3}
\end{align}
where the fine grid length scale $\Delta x$ would be a fraction of $h$. To compute derivatives at boundaries left/right stencils\cite{fornberg1988generation} can be exploited. This suggests primal and dual grids used in Hermite interpolation can be coarser than the original grid, where data is available, without loss of information. In contrast, this would not be the case with tricubic spline interpolation.
Note that the expressions \eqref{eq:fd1}-\eqref{eq:fd3} are fourth order finite difference approximations for the first and the second derivatives.
Alternatively, projection \cite{beznosov2021hermite} or interpolation methods could be used to find the coefficients in \eqref{coef1D}.

\subsection{\label{sssec:oneD} Hermite interpolation in 1D}
\newcommand{\xr}{x_{\rm R}}
\newcommand{\xl}{x_{\rm L}}
Next, we consider an interval bounded by two neighboring grid points $x_i$, $x_{i+1} \in \Omega_p$ and construct the unique local Hermite interpolating polynomial of degree $(2m+1)$, centered on the dual grid point $x_j \in \Omega_d, $ $j = i+\frac 12$. The interpolating polynomials can be written in Taylor form
\begin{equation}
u_{j}(x)  = \sum_{l = 0}^{2m+1} \hat u_{j,l} \left( \frac{x - x_{j}}{h}\right)^l, \ \ j = \frac 12, \dots, N-\frac 12. \label{hcoef1D}
\end{equation}
The interpolating polynomial $u_{i+\frac 12}$ is determined by the local interpolation conditions
\begin{equation}
  \left(\frac{{\rm d}^l u_{j}}{{\rm d}x^l} = \frac{{\rm d}^l u_{i}}{{\rm d}x^l}\right)\bigg|_{x = x_i}, \ \
  \left(\frac{{\rm d}^l u_{j}}{{\rm d}x^l} = \frac{{\rm d}^l u_{i+1}}{{\rm d}x^l}\right)\bigg|_{x = x_{i+1}}, \ \ l = 0,\dots, m, \ \  j = i + \frac 12,\label{eq:int_cond}
\end{equation}
and we find the coefficients in \eqref{hcoef1D} by forming a generalized Newton table as described in \cite{hagstrom2015solving}.

\subsection{\label{sssec:Differentiation_and_integration} Differentiation and integration}
The piece-wise Hermite interpolation polynomials have a local monomial
basis and thus are easy to differentiate and integrate analytically by
shifting and rescaling the coefficients. The resulting derivatives
can be expressed in Taylor form as well. The derivative of the polynomial
\eqref{hcoef1D} is
\begin{equation}
\frac {{\rm d}^n}{dx^n} u_{j}(x) = \sum_{l = 0}^{2m+1 - n} \frac {(l+1) \dots (l+n)}{h^n}\hat u_{j, l+n} \left(\frac {x - x_j} h\right)^l. \label{eq:diff}
\end{equation}
\begin{comment}
  Its antiderivative, written as the integral from the left endpoint $x_i$ to a variable point $x \in [x_i, x_{i+1}]$ can be expressed in the same (monomial) basis as
  \begin{equation}
  \int_{x_i}^x u_{i+\frac 12}(x')\ {\rm d}x' = s
  \sum_{l = 1}^{2m+2} \frac {h}{l}\hat u_{i+\frac 12, l-1} \left(-\frac 12\right)^l
  +
  \sum_{l = 1}^{2m+2} \frac {h}{l}\hat u_{i+\frac 12, l-1} \left(\frac {x - x_{i + \frac 12}} h\right)^l. \label{eq:int_l}
  \end{equation}
  Similarly, the integral from a variable point $x \in [x_i, x_{i+1}]$ to the right endpoint $x_{i+1}$ is
  \begin{equation}
  \int_{x}^{x_{i+1}} u_{i+\frac 12}(x')\ {\rm d}x' =
  \sum_{l = 1}^{2m+2} \frac {h}{l}\hat u_{i+\frac 12, l-1} \left(\frac 12\right)^l
  -
  \sum_{l = 1}^{2m+2} \frac {h}{l}\hat u_{i+\frac 12, l-1} \left(\frac {x - x_{i + \frac 12}} h\right)^l. \label{eq:int_r}
  \end{equation}
\end{comment}
Further, its indefinite integral is
\begin{equation}
\int u_{j}(x)\ dx = \sum_{l = 1}^{2m+2} \frac {h}{l}\hat u_{j, l-1} \left(\frac {x - x_{j}} h\right)^l. \label{eq:int_r}
\end{equation}
and integrating over dual cell endpoints, we obtain
\begin{equation}
\int_{x_j - \frac h2}^{x_j + \frac h2} u_{j}(x)\ {\rm d}x =
\sum_{k = 0}^{m} \frac {h}{2k+1}\hat u_{j, 2k} \left(\frac 12\right)^{2k}. \label{eq:int_cell}
\end{equation}
Note that the expressions \eqref{eq:diff}-\eqref{eq:int_cell} are exact. By combining \eqref{eq:int_r} and \eqref{eq:int_cell}, we can compute the approximation to the integral of $u$ from an arbitrary dual grid point $x_{j}$ to an arbitrary point $y \in [x_k - \frac h2, x_k + \frac h2]$ with $j < k$ as
\begin{widetext}
\begin{align}
\int_{x_j}^x  f(x')\ {\rm d}x' &\approx \int_{x_j}^{x_j+\frac h2} u_{j}(x)\ {\rm d}x
+
\sum_{j' = j+1}^{k-1} \int_{x_{j'} - \frac h2}^{x_{j'} + \frac h2} u_{j'}(x)\ {\rm d}x +
\int_{x_k - \frac h2}^x u_{j}(x)\ {\rm d}x \nonumber \\
&=
\sum_{l = 1}^{2m+2} \frac hl \left[
   \sum_{j' = j}^{k-1}
    \hat u_{j', l-1}\left( \frac 12 \right)^l -
   \sum_{j' = j+1}^{k}
   \hat u_{j', l-1}
\left(-\frac 12 \right)^l\right] +
\sum_{l = 1}^{2m+2} \frac {h}{l}\hat u_{k, l-1} \left(\frac {x - x_{k}} h\right)^l. \label{eq:int_full}
\end{align}
\end{widetext}
This defines an approximation to the anti-derivative of $u$ written in Taylor form.
\subsection{\label{ssec:Higher_dimensions}Higher dimensions}

In higher dimensions the approximations to $u$ take the form of centered tensor product Taylor polynomials. In the case of cylindrical coordinates -- which are often considered in tokamak fusion reactor descriptions -- the coefficients would be of the form $\hat u^{j_1,j_2}_{l,s;k}$, with superscripts representing cell indices of the $(R, Z)$ poloidal plane. Further, the first two subscripts represent powers
in the two spatial directions of the poloidal plane, and the third represents the angular (toroidal) coordinate $\varphi$. Since the $\varphi$-direction is periodic, a smooth, easy to differentiate approximation can be achieved by using trigonometric polynomials. In particular, a single, constant coefficient in the $\varphi$ direction would correspond to axisymmetric fields.

Altogether, we consider tensor product Taylor-trigonometric polynomials given by
\begin{equation}
  u^{i_1, i_2} (R,\varphi,Z) = \sum_{k = 0}^{N_\varphi}\sum_{l = 0}^{m_R}\sum_{s=0}^{m_Z} \hat u^{i_1, i_2}_{l,s;k} \left(\frac{R - R_{i_1}}{h_R}\right)^l\left(\frac{Z - Z_{i_2}}{h_Z}\right)^s T_k(\varphi),\ \   i_1 = 0,\dots,N_R,\ \ i_2 = 0, \dots, N_Z, \label{Taylor-trig}
\end{equation}
where $T_k$ is a desired trigonometric polynomial basis. In this work we use a real discrete Fourier transform of the field data in $\varphi$, with an associated basis
\begin{equation}
  T_k(\varphi) =
  \begin{cases}
    \cos(k\varphi), & k = 0,\dots, \frac{N_\varphi}2,\\
    \sin\left(\left(k - \frac{N_\varphi}2\right)\varphi \right), & k = \frac{N_\varphi}2+1,\dots, {N_\varphi}.\\
  \end{cases} \label{eq:trig_basis}
\end{equation}
A tensor product version of the finite differences \eqref{eq:fd1}-\eqref{eq:fd3} can be used to determine the poloidal coefficients. Further, we compute toroidal coefficients working out the finite differences in several equidistant poloidal planes and performing a Fourier transform for each index set $(i_1,i_2,l,s)$.

To perform Hermite interpolation in the two-dimensional poloidal plane, we apply sequences of one-dimensional Hermite interpolations as described in Section~(\ref{sssec:oneD}) in two directions as follows. First, we form the polynomials $p_{j_1, i_2}(R,Z)$ centered at each cell edge by performing $(m_Z+1)$ $R-$interpolations, satisfying local interpolation conditions \eqref{eq:int_cond} for $R$, so that
\begin{align}
  &\left[\frac {{\rm d}^n}{{\rm d} R^n}\sum_{l = 0}^{2m_R+1}\hat u^{{j_1},i_2}_{l,s;k} \left(\frac{R - R^d_{j_1}}{h_R}\right)^l\right]\Bigg|_{R = R^p_{i_1}} =
  \left[\frac {{\rm d}^n}{{\rm d} R^n}\sum_{l = 0}^{m_R}\hat u^{{i_1},i_2}_{l,s;k} \left(\frac{R - R^p_{i_1}}{h_R}\right)^l\right]\Bigg|_{R = R^p_{i_1}}, \label{eq:int_cond1}\\
  &\left[\frac {{\rm d}^n}{{\rm d} R^n}\sum_{l = 0}^{2m_R+1}\hat u^{{j_1},i_2}_{l,s;k} \left(\frac{R - R^d_{j_1}}{h_R}\right)^l\right]\Bigg|_{R = R^p_{i_1+1}} =
  \left[\frac {{\rm d}^n}{{\rm d} R^n}\sum_{l = 0}^{m_R}\hat u^{{i_1+1},i_2}_{l,s;k} \left(\frac{R - R^p_{i_1+1}}{h_R}\right)^l\right]\Bigg|_{R = R^p_{i_1+1}}, j_1 = i_1 + \frac 12. \label{eq:int_cond2}
\end{align}
This results in $(2m_R+1) \times (m_Z+1)$-degree polynomials centered on intermediate grid points $(R_{j_1}, Z_{i_2})$ located on vertical cell edges. Superscripts $d$ and $p$ indicate dual and primal grid points, respectively.

We then form polynomials $p_{j_1, j_2}(R,Z)$ centered on dual grid points in cell centers by performing an additional $(2m_R+2)$ $Z-$interpolations per dual grid point satisfying \eqref{eq:int_cond1}, \eqref{eq:int_cond2}. The above operations can be performed for each $\varphi$-mode independently and in parallel, and the resulting piece-wise polynomial becomes
\begin{equation}
  u^{j_1,j_2} (R,\varphi,Z) = \sum_{k = 0}^{N_\varphi}\sum_{l = 0}^{2m_R+1}\sum_{s=0}^{2m_Z+1} \hat u^{j_1,j_2}_{l,s;k} \left(\frac{R - R^d_{j_1}}{h_R}\right)^l\left(\frac{Z - Z^d_{j_2}}{h_Z}\right)^s T_k(\varphi),\ \   j_1 = \frac 12,\dots,N_R - \frac 12,\ \ j_2 = \frac 12, \dots, N_Z-\frac 12. \label{eq:2dpoly}
\end{equation}

Differentiation and integration of \eqref{eq:2dpoly} with respect to $R$ and $Z$ can be performed analytically, as described in Section~(\ref{sssec:Differentiation_and_integration}). In this work we also need differentiation in $\varphi$, which can be trivially performed within the trigonometric basis \eqref{eq:trig_basis}.

\section{\label{sec:Flux_reconstruction} Vector potential reconstruction}
In the context of a right-hand side containing a discrete magnetic field $\mathbf{B}$, we interpolate $\mathbf{B}$ indirectly through its vector potential $\mathbf{A}$, in order to ensure a divergence-free construction. For this purpose, we first recall that in cylindrical coordinates $(R,\varphi,Z),$ the curl of the vector potential $\mathbf{A}$ gives rise to the magnetic field in component form:
\begin{align}
  B_R & = \frac{1}{R}\frac{\partial A_Z}{\partial \varphi} - \frac{\partial A_\varphi}{\partial Z} \\
  B_\varphi & = \frac{\partial A_R}{\partial Z} - \frac{\partial A_Z}{\partial R} \\
  B_Z & = \frac{1}{R}\frac{\partial}{\partial R}\left(RA_\varphi\right) - \frac{1}{R}\frac{\partial A_R}{\partial\varphi}.
\end{align}
\subsection{\label{ssec:rat_poly_approximation} Rational approximation}
In the case of an axisymmetric magnetic field, we have
\begin{align}
  \psi(R,Z) & = RA_\varphi \\
  G(R,Z) & = RB_\varphi,
\end{align}
where both $\psi(R,Z)$ and $G(R,Z)$ are well approximated by polynomials
of $R$ and $Z.$ Because of this, $A_\varphi = \psi(R,Z)/R$ and
$B_\varphi = G(R,Z)/R$ are naturally expressed as rational
polynomials. This implies that we should not interpolate $A_\varphi$
and $B_\varphi$ directly, and instead interpolate the quantities $RA_\varphi$ and $RB_\varphi.$
Since the magnetic field's axisymmetric $B_\varphi$-component dominates the overall
magnetic field in tokamaks, the last statement can also be applied to the full quantities
$RA_\varphi(R,\varphi,Z)$ and $RB_\varphi(R,\varphi,Z)$, possibly including
non-axisymmetric components. For a general three dimensional toroidal
field, both $\psi(R,\varphi,Z)$ and $G(R,\varphi,Z)$ will depend on
$\varphi$,
\begin{align}
  \psi(R,\varphi,Z) & = RA_\varphi(R,\varphi,Z) \\
  G(R,\varphi,Z) & = RB_\varphi(R,\varphi,Z).
\end{align}
Given $\psi$ and $G$ -- which now simply correspond to rescaled versions of $A_\varphi$ and $B_\varphi$ rather than the original poloidal and toroidal flux functions in the context of axisymmetric magnetic fields -- the relation $\mathbf{B}=\nabla\times\mathbf{A}$ can be formulated as
\begin{align}
  B_R & = \frac{1}{R}\frac{\partial A_Z}{\partial \varphi} - \frac{1}{R}\frac{\partial \psi }{\partial Z} \\
  G & = \frac{\partial R A_R}{\partial Z} - R \frac{\partial A_Z}{\partial R} \label{eqn:GfromA}\\
  B_Z & = \frac{1}{R}\frac{\partial\psi}{\partial R} - \frac{1}{R}\frac{\partial A_R}{\partial\varphi}.
\end{align}
Gauge freedom allows us to set one of $\mathbf{A}$'s components to vanish, and the form of the $G$ equation \eqref{eqn:GfromA} implies a good candidate to be
\begin{align}
A_Z = 0,
\end{align}
together with a new variable
\begin{align}
\chi(R,\varphi,Z) \equiv R A_R(R,\varphi,Z).
\end{align}
From \eqref{eqn:GfromA}, it is easy to see that $\chi(R,\varphi,Z)$ is well approximated by a polynomial in $R$ and $Z,$
because $G$ is.
Finally, the reconstruction of $\mathbf{A},$ which is in the form of rational poloynomials
\begin{equation}
  A_R = \frac{\chi}{R}, \quad
  A_\varphi = \frac{\psi}{R},\quad
  A_Z = 0, \label{eq:vp}
\end{equation}
is done by integrating
\begin{align}
  - \frac{1}{R}\frac{\partial \psi }{\partial Z} & = B_R \\
  \frac{\partial \chi}{\partial Z} & = G = R B_\varphi\\
  \frac{1}{R}\frac{\partial\psi}{\partial R} - \frac{1}{R^2}\frac{\partial \chi}{\partial\varphi} & = B_Z,
\end{align}
that is we compute $\psi$ and $\chi$ through the magnetic field $(B_R, \; G\!=\!R B_\varphi, \; B_Z)$ as defined on its given discrete grid. Integrating the first two equations with respect to $Z$, we obtain
\begin{align}
  \psi(R,\varphi,Z) &=  - \int_{Z_0}^Z R B_R(R, \varphi, Z')\ {\rm d}Z' + C_1(R,\varphi),\\
  \chi(R,\varphi,Z) &=  \int_{Z_0}^Z R B_\varphi(R, \varphi, Z')\ {\rm d}Z' + C_2(R,\varphi).
  \end{align}
  Substituting these into the last component equation for $B_Z$ and using $\nabla\cdot\mathbf{B}=0,$  one finds
  \begin{align}
    \frac{1}{R}\frac{\partial\psi}{\partial R}
    - \frac{1}{R^2}\frac{\partial\chi}{\partial \varphi} - B_Z
    = & -  \int_{Z_0}^Z \frac{1}{R}\frac{\partial}{\partial R} \left[R B_R(R, \varphi, Z')\right] {\rm d}Z'
    + \frac{1}{R}\frac{\partial C_1}{\partial R} \nonumber \\
    & - \frac{1}{R} \int_{Z_0}^Z \frac{\partial}{\partial \varphi} B_\varphi(R, \varphi, Z')\ {\rm d}Z'
    - \frac{1}{R^2}\frac{\partial C_2}{\partial\varphi} - B_Z \nonumber \\
    = & \int_{Z_0}^Z \frac{\partial B_Z(R,\varphi,Z')}{\partial Z'}\ {\rm d}Z' + \frac{1}{R}\frac{\partial C_1}{\partial R}
    - \frac{1}{R^2}\frac{\partial C_2}{\partial\varphi} - B_Z \nonumber \\
    = & B_Z(R,\varphi,Z) - B_Z(R,\varphi,Z_0) + \frac{1}{R}\frac{\partial C_1}{\partial R}
    - \frac{1}{R^2}\frac{\partial C_2}{\partial\varphi} - B_Z(R,\varphi,Z) \nonumber \\
    = & - B_Z(R,\varphi,Z_0) + \frac{1}{R}\frac{\partial C_1}{\partial R}
    - \frac{1}{R^2} \frac{\partial C_2}{\partial\varphi} \nonumber \\
    = & 0.
  \end{align}
  This equation can be satisfied by simply taking $C_2=0$ and finding $C_1$ from
  \begin{align}
  C_1(R,\varphi) = \int_{R_0}^R R' B_Z(R',\varphi,Z_0) \ {\rm d}R'.
  \end{align}

Finally, gathering all pieces together, we have
\begin{align}
\psi(R,\varphi,Z) &= - \int_{Z_0}^Z R B_R(R, \varphi, Z') {\rm d}Z' + \int_{R_0}^R R' B_Z(R',\varphi,Z_0) \ {\rm d}R', \label{rec:psi}\\
\chi(R,\varphi,Z) &=  \int_{Z_0}^Z R B_\varphi(R, \varphi, Z') {\rm d}Z'. \label{rec:chi}
\end{align}

\subsection{\label{ssec:Reconstruction} Reconstruction procedure}
To evaluate \eqref{rec:psi} and \eqref{rec:chi}, we consider a single poloidal cross-section of a torus representing a
tokamak fusion reactor, discretized by a rectangular grid $\Omega_p = (R_{i_1},
Z_{i_2})$ with points
\begin{align}
  &R_{i_1} = R_{\min} + i_1h_R,  \label{eq:fd21}\\
  &Z_{i_2} = Z_{\min} + i_2h_Z,  \label{eq:fd22}
\end{align}
and employ the finite differences \eqref{eq:fd1}-\eqref{eq:fd3} on a finer grid to compute fields and derivatives, with finite difference stencils formed around $\Omega_p$ grid points. We form Taylor polynomials $(u^{i_1,i_2}$,
$v^{i_1,i_2}$, $w^{i_1,i_2})$ on the primal grid $\Omega_p$, which approximate the possibly non-divergence free discrete field components, $B_R$, $B_\varphi$ and $B_Z$, respectively. The polynomials are of the form \eqref{Taylor-trig}. Next, using Section~(\ref{ssec:Higher_dimensions}), we obtain interpolation polynomials $(u^{j_1,j_2}$, $v^{j_1,j_2}$, $w^{j_1,j_2})$.
This will give a continuous representation of the components of $\mathbf B$, which has $m_R$ and $m_Z$ derivatives with guaranteed continuity both locally and across cell interfaces.

The resulting piece-wise polynomial representation can be integrated analytically, see Section~(\ref{sssec:Differentiation_and_integration}), thereby providing a path to compute
$\psi$ and $\chi$ in \eqref{rec:psi} and \eqref{rec:chi}. Applying \eqref{eq:int_full} separately in the $R$ and $Z$ directions, we obtain
\begin{align}
  \psi(R, \varphi, Z) &\approx p^{j_1,j_2}(R,\varphi,Z) = -\int_{Z_c}^Z u(R,\varphi, Z')\ {\rm d} Z' +
    \int_{R_c}^R w(R',\varphi, Z_c)\ {\rm d} R', \label{eq:psi_approx}\\
  \chi(R, \varphi, Z) &\approx q^{j_1,j_2}(R,\varphi,Z) = \int_{Z_c}^Z v(R,\varphi, Z')\ {\rm d} Z', \quad (R,Z) \in [R^d_{j_1}- \frac {h_R}2, R^d_{j_1} + \frac {h_R}2] \times [Z^d_{j_2} - \frac {h_Z}2, Z^d_{j_2} + \frac {h_Z}2], \label{eq:chi_approx}
\end{align}
where $p^{j_1, j_2}$ and $q^{j_1, j_2}$ are $(2m_R+2)\times(2m_Z+2)\times (N_\varphi+1)$ and $(2m_R+1)\times(2m_Z+2) \times(N_\varphi+1)$ piece-wise  Taylor-trigonometric polynomials defined locally on each cell.

\subsection{\label{ssec:Field_evaluation}Field evaluation}
\begin{figure}
  \centering
\includegraphics[width=0.24\columnwidth]{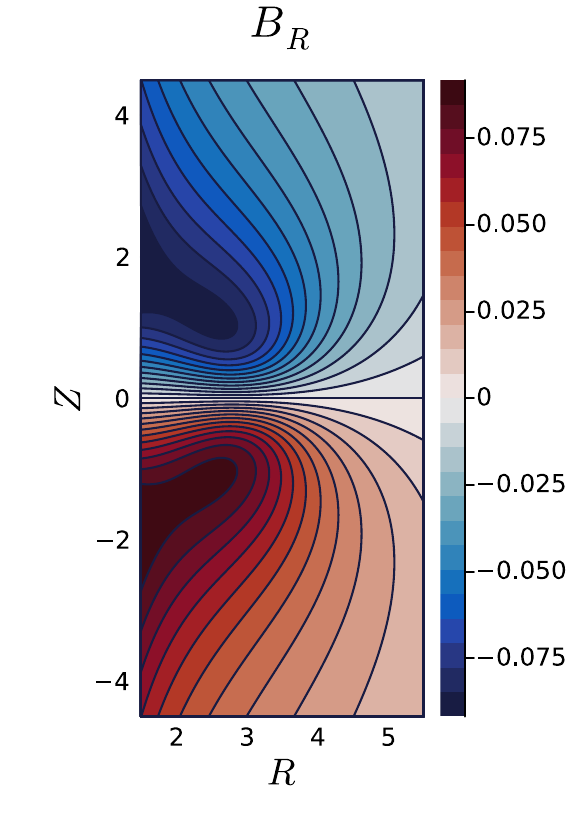}
\includegraphics[width=0.24\columnwidth]{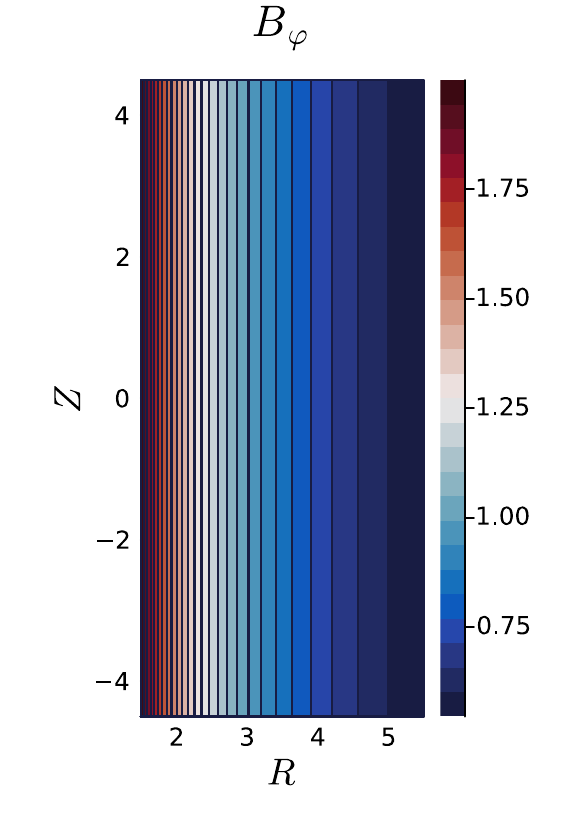}
\includegraphics[width=0.24\columnwidth]{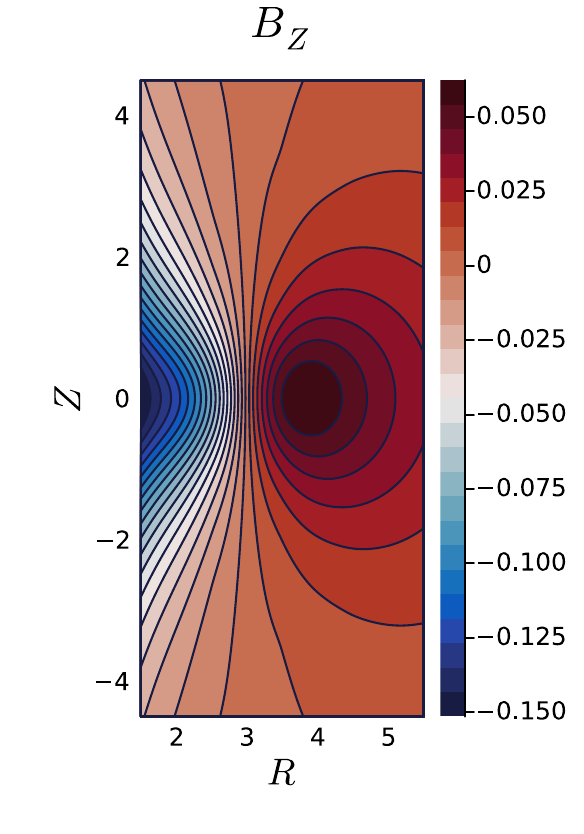}
\includegraphics[width=0.24\columnwidth]{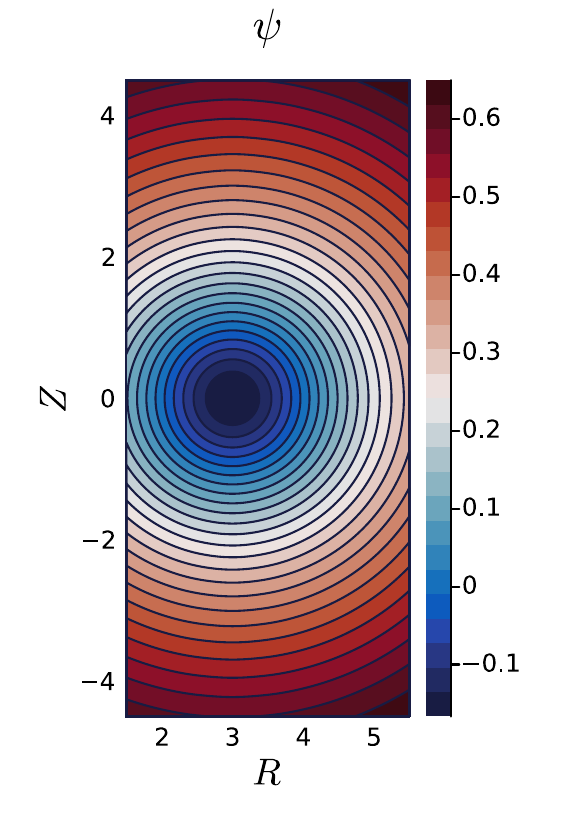}
\caption{\label{fig:fields} 9th order Hermite interpolation of analytically constructed magnetic field \eqref{eq:B_from_q}, and reconstructed poloidal flux function on a coarsest mesh. The number of grid points in each direction is given by $N_R = 5$,  $N_Z = 10$, and the resulting field approximation is continuous, differentiable and divergence free.
}
\end{figure}
Finally, given the approximated vector potential \eqref{eq:vp} through \eqref{eq:psi_approx} and \eqref{eq:chi_approx}, we can apply a curl using \ref{sssec:Differentiation_and_integration} to obtain an interpolated expression for the magnetic field, given by
\begin{align}
  &B_R \approx \frac {u(R, \varphi, Z)} R, \quad
  B_\varphi \approx \frac {v(R, \varphi, Z)} R, \label{eq:int_expr1} \\
  &B_Z \approx \frac 1R \left(-\int_{Z_c}^Z\frac{\partial u} {\partial R}(R, \varphi, Z')\  {\rm d}Z'
  + w(R, \varphi, Z_c) \right)
  -\frac 1{R^2} \int_{Z_c}^Z\frac{\partial v} {\partial \varphi}(R, \varphi, Z')\ {\rm d}Z' \label{eq:int_expr2},
\end{align}
where
\begin{align}
&u(R, \phi, Z) = u^{j_1, j_2}(R, \phi, Z),\ \ 
v(R, \phi, Z) = v^{j_1, j_2}(R, \phi, Z),\ \
w(R, \phi, Z) = w^{j_1, j_2}(R, \phi, Z),\\ 
& R_{j_1} - \frac {h_R}2 < R < R_{j_1} + \frac {h_R} 2, \\ 
& Z_{j_2} - \frac {h_Z}2 < Z < Z_{j_2} + \frac {h_R} 2, \ \ (R,Z) \in \Omega_d.
\end{align}
Derivatives can be taken by using \eqref{eq:diff} and trigonometric identities. Since all differentiation is exact, the resulting field's divergence is $0$. Moreover if $m_R > 1$ and $m_Z > 1$, its gradient and curl  are continuous.

\section{\label{sec:Numerical_experiments}Numerical experiments}
\subsection{\label{ssec:Convergence} Analytical circular flux surfaces: rates of convergence}
Consider a axisymmetric magnetic field represented as
\begin{equation}
  \mathbf{B} = R B_\varphi \nabla \varphi + \nabla \varphi \times \nabla \psi, \ \ RB_\varphi = \text{const},
\end{equation}
where here we set $RB_\varphi$ equal to a constant for simplicity, rather
than the usual flux function considered for axisymmetric field configurations.
The magnetic field's components are defined with respect to the magnetic axis at \((R,Z)
= (3,0)\), where the length has been normalized by a minor radius $a=2$~meters,
so that the aspect ratio of \(3\) corresponds to the value of ITER. We then set
\begin{equation}
  RB_R (R, Z) = -\frac {Z}{q(R,Z)},\quad
  RB_\varphi (R, Z) = 3, \quad
  RB_Z (R, Z) = \frac {R - 3}{q(R,Z)}, \label{eq:B_from_q}
\end{equation}
and for the results presented in this section, the safety factor's profile is given by
\begin{equation}
  q(R,Z) = q_0 + q_2 (R - 3)^2 + q_2Z^2,\quad q_0 = 2.0, \quad q_2 = 2.1.
\end{equation}
The poloidal flux function is $\psi(R,Z) = \frac{\log(q)}{2q_2}$.

\begin{figure}
  \centering
\includegraphics[width=0.49\columnwidth]{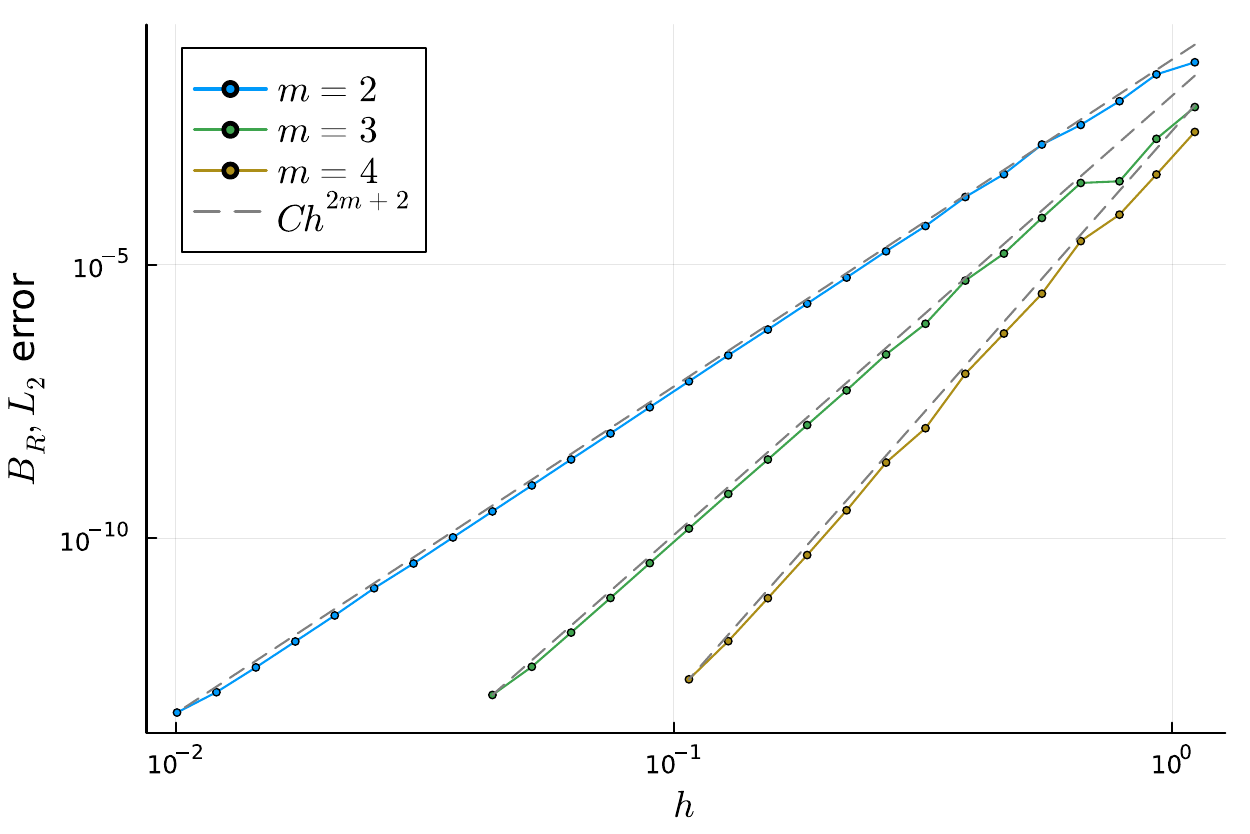}
\includegraphics[width=0.49\columnwidth]{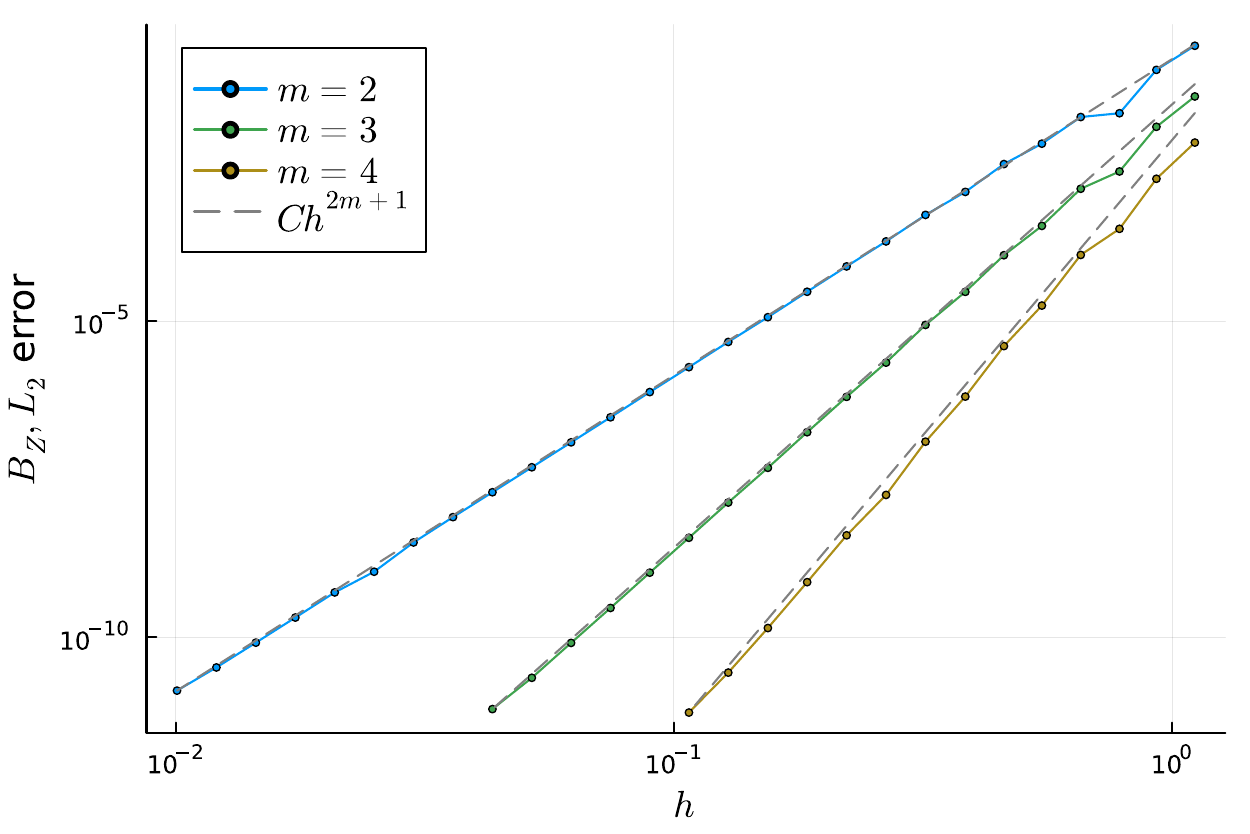}
\includegraphics[width=0.49\columnwidth]{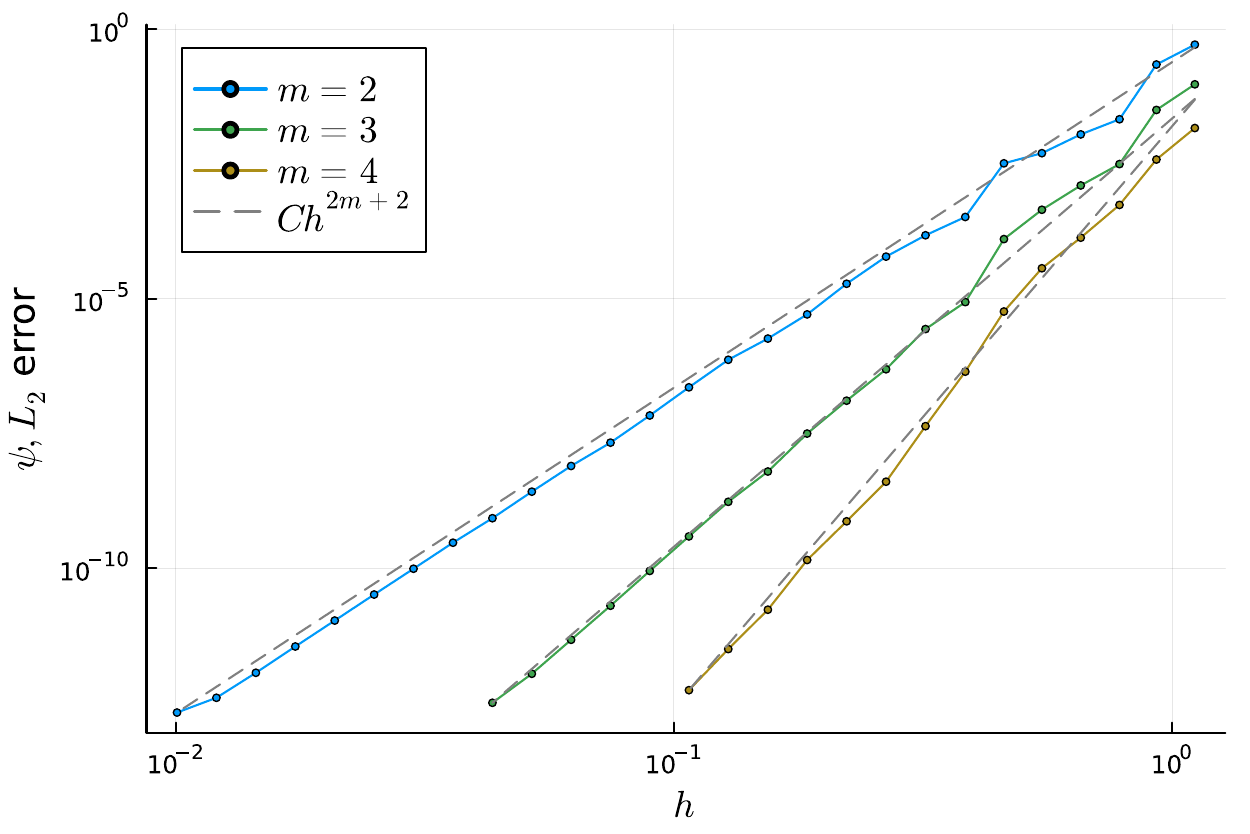}
\includegraphics[width=0.49\columnwidth]{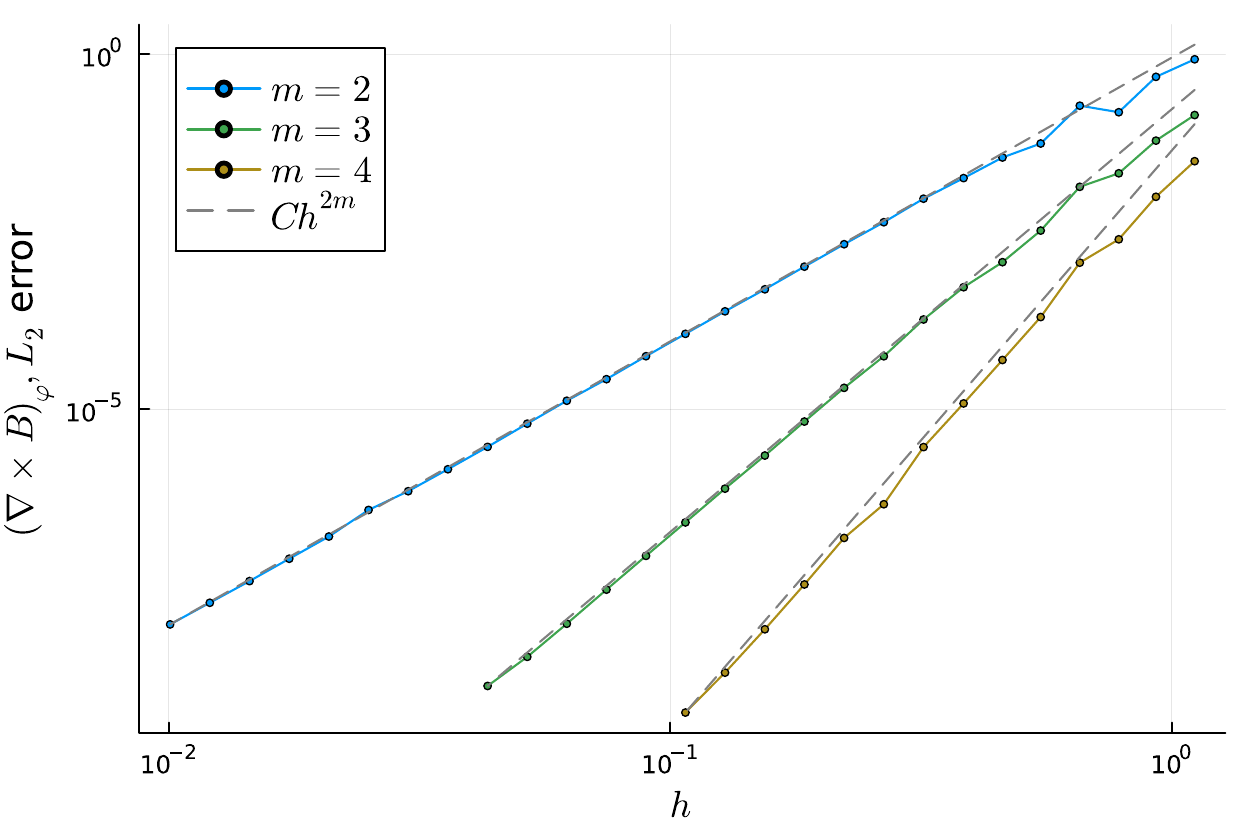}
\caption{\label{fig:convergence} Hermite interpolation convergence
  for analytic magnetic fields. Solid lines display the relative
  error arising when computing magnetic field components and poloidal flux
  function as a function of grid spacing used in the Hermite
  interpolation. Dashed lines display the expected convergence
  rates. $\psi$ (lower left) and $B_R$ (upper left) approximations
  converge at 6th, 8th and 10th order. $B_\varphi$ (not shown) does
  not have an error because $RB_\varphi$ is a constant and thus is
  exactly represented in the polynomial basis. $B_Z$ (upper right)
  converges at a reduced order (5th, 7th and 9th). The reason for
  this is that $B_Z$ is approximated by a derivative of $u$
  with respect to $R$, which truncates at the $2m_R = 4,6,8$
  coefficients. Further, $\nabla B$ (lower right) converges at 4,6
  and 8th order.  }
\end{figure}
We consider a rectangular domain $[1,6] \times [-5, 5]$, which is discretized with a $n_R\times n_Z$ regular grid. The number of grid points of the Hermite grid is set to $N = (n_R-1)/4$ and $M = (n_Z-1)/4$ in the $R$- and $Z$-directions, respectively. Starting with $n_R = 21$ and $n_Z = 41$ grid points in each direction, we measure the error and increase the resolution while keeping the same aspect ratio.
For the error measurement, we evaluate an approximation on a finer grid oversampled with 10000 grid points per unit square and compute the $l_2$ norm of the difference to an exact solution evaluated on the same mesh, relative to the $l_2$-norm of the solution. The convergence plot is shown in FIG.~(\ref{fig:convergence}). Solid lines display the relative errors appearing in the magnetic field component- and poloidal flux function computations, as functions of grid spacing $h_Z$ in the Hermite interpolation. As indicated by the plot, the $\psi$ and $B_R$ approximations converge at fifth order. $B_Z$ converges at an order decreased by 1. The reason for this decrease is that $B_Z$ is approximated using a derivative of $u$ with respect to $R$, which truncates at the $2m_R$ coefficient.

\subsection{\label{ssec:First_order_mimetic_finite_difference_fields} Conservation in guiding center motion}
We consider the relativistic guiding center equations\cite{mcdevitt2019avalanche,cary2009hamiltonian} in the coordinate system $(p, \xi, \mathbf{X})$, in the absence of collisions, radiation and the electric field. In terms of dimensionless variables, they can be written as
\begin{align}
     &\frac {{\rm d}p}{dt} = 0, \\
     &\frac {{\rm d}\xi}{dt} = - \frac {p(1 - \xi^2)} {2\gamma \mathbf{B}^\ast_\parallel} \mathbf{B}^\ast \cdot \nabla \ln {B}, \\
     &\frac {{\rm d}\mathbf{X}}{{\rm d}t} = \frac {\xi p}{\gamma} \frac{\mathbf{B}^\ast}{\mathbf{B}^\ast_\parallel}
   - \frac 1{\omega_c} \frac{p^2 (1-\xi^2)}{2\gamma \mathbf{B}_\parallel^\ast} \vec {\mathbf{b}} \times \nabla \ln {B},
\end{align}
where
\begin{align}
     \gamma = \sqrt{1+p^2}, \quad \mathbf{B}^\ast = \mathbf{B} - \xi p \nabla \times \vec {\mathbf{b}}, \quad
     \mathbf{B}^\ast_\parallel = \mathbf{B}^\ast \cdot \vec{\mathbf{b}}.
\end{align}

Here time $t$ is normalized by the tokamak's minor radius over the speed of light, $a/c$. $p$ is the particle momentum normalized by $mc$, $\xi = \frac {p_\parallel} p$ is the pitch, $\mathbf{X}$ is the guiding center position, $\vec {\mathbf{b}} = \mathbf{B} / B$ is the direction of magnetic field, and $\omega_c$ is the cyclotron frequency of a charged particle positioned at the magnetic axis.
\begin{figure}
  \centering
\includegraphics[width=0.7\columnwidth]{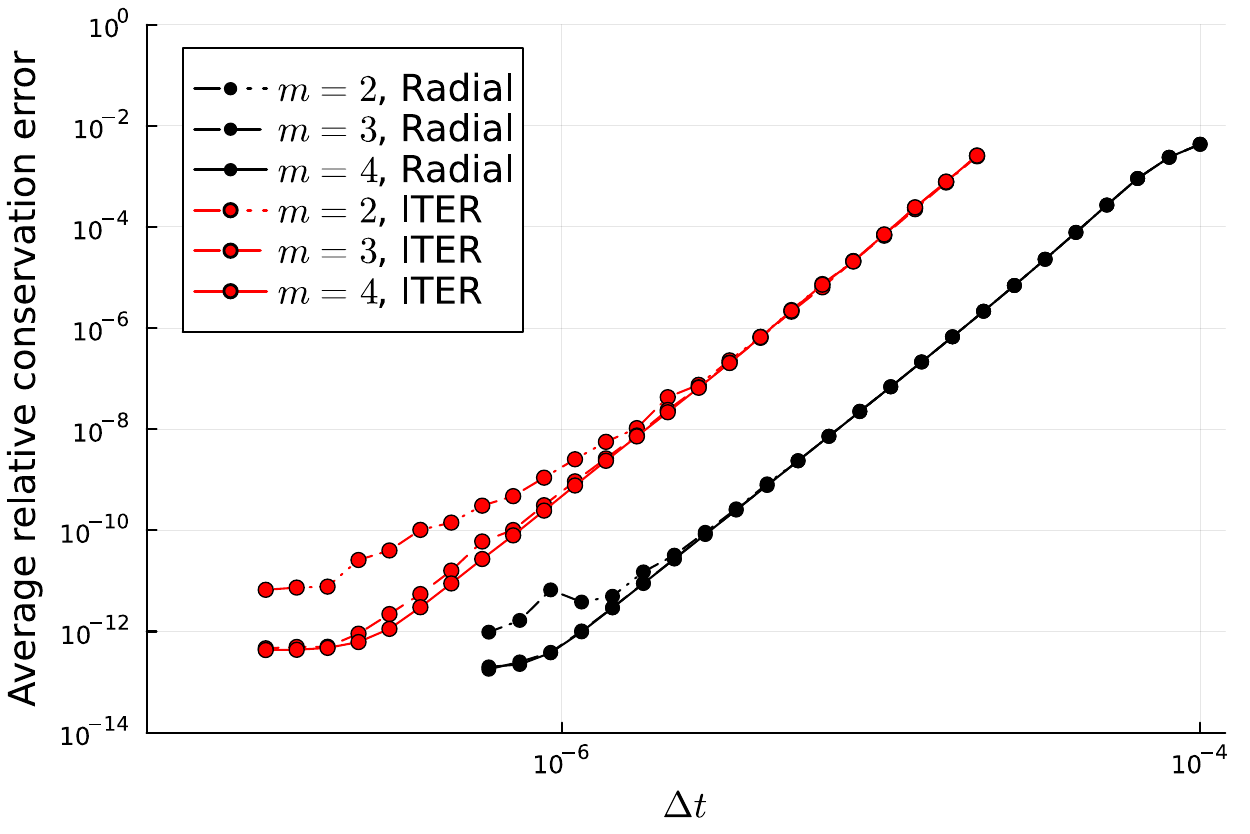}
\caption{\label{fig:cons_comp} Guiding center equations conservation-convergence test. Average, relative conservation error of $p_\varphi$ and $\mu$ as a function of the timestep. The timestep is normalized by 0.5 ms, which is the simulation's final time. Lines display the conservation errors of the 5th order Runge-Kutta method from the Dormand-Prince embedded pair of solutions. Red lines display the conservation error when background fields are set to the ITER-type VDE simulation, as displayed in FIG.~(\ref{fig:fields_mfd}). Black lines display the conservation error when background fields are set to analytic fields with radial flux surfaces, as displayed in FIG.~(\ref{fig:fields}). Different line styles represent the number of derivatives interpolated using the Hermite method for field reconstruction.}
\end{figure}
\begin{figure}
  \centering
\includegraphics[width=0.19\columnwidth, trim= {5pt 5pt 15pt 10pt}]{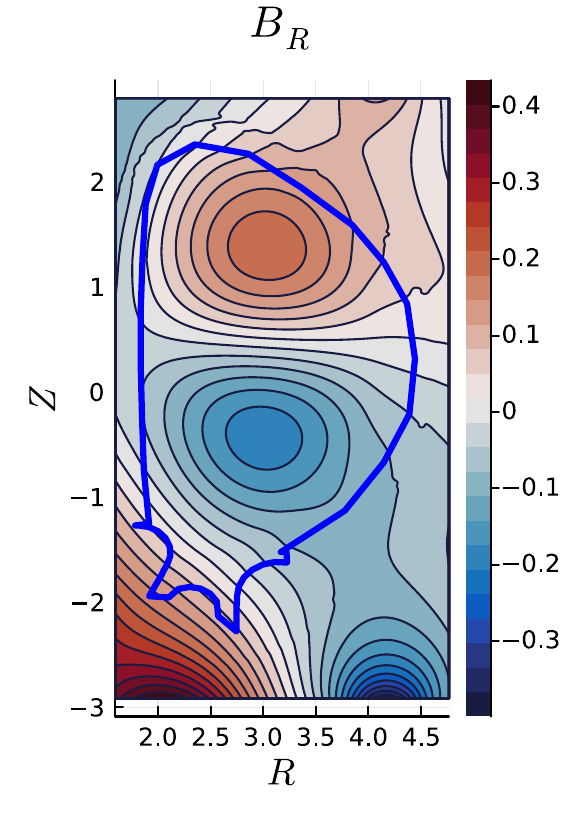}
\includegraphics[width=0.19\columnwidth, trim= {5pt 5pt 15pt 10pt}]{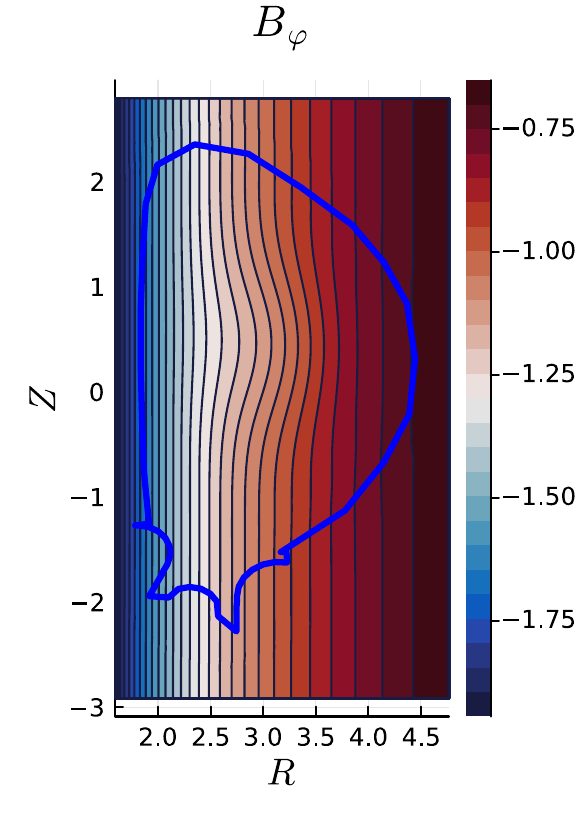}
\includegraphics[width=0.19\columnwidth, trim= {5pt 5pt 15pt 10pt}]{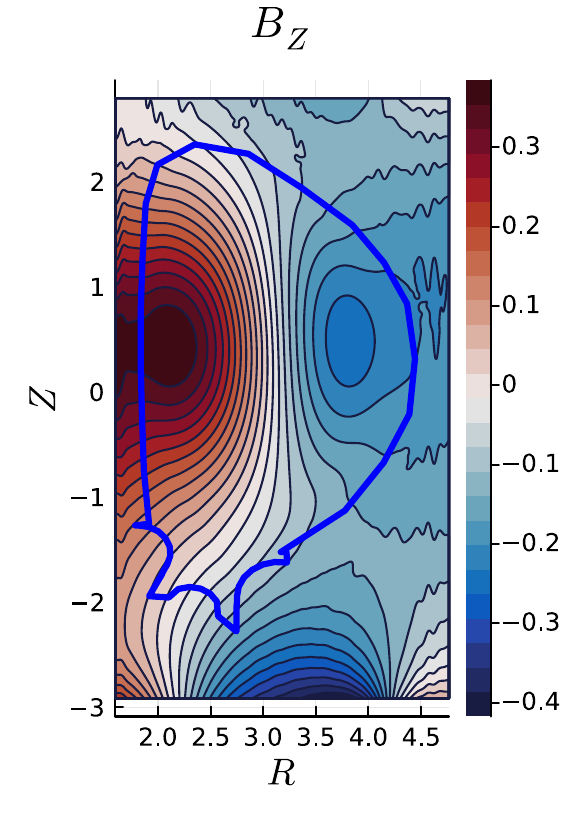}
\includegraphics[width=0.19\columnwidth, trim= {5pt 5pt 15pt 10pt}]{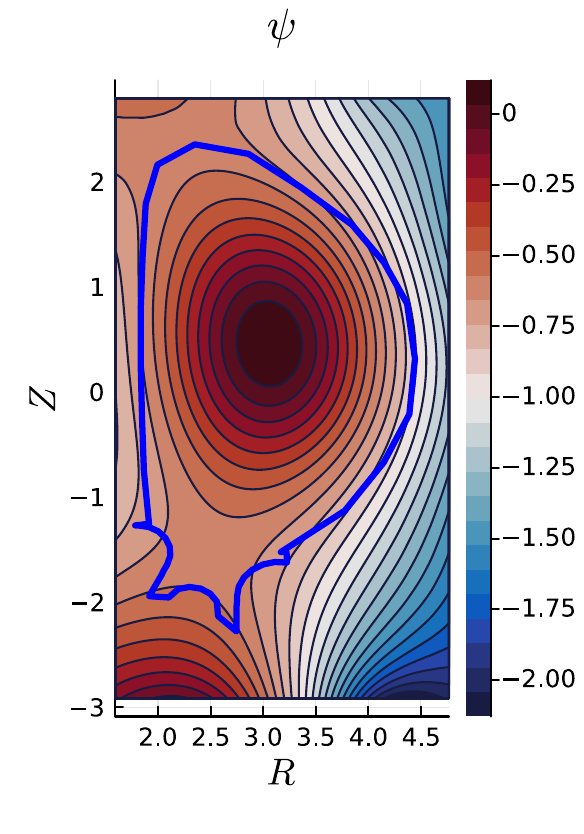}
\includegraphics[width=0.19\columnwidth, trim= {5pt 5pt 15pt 10pt}]{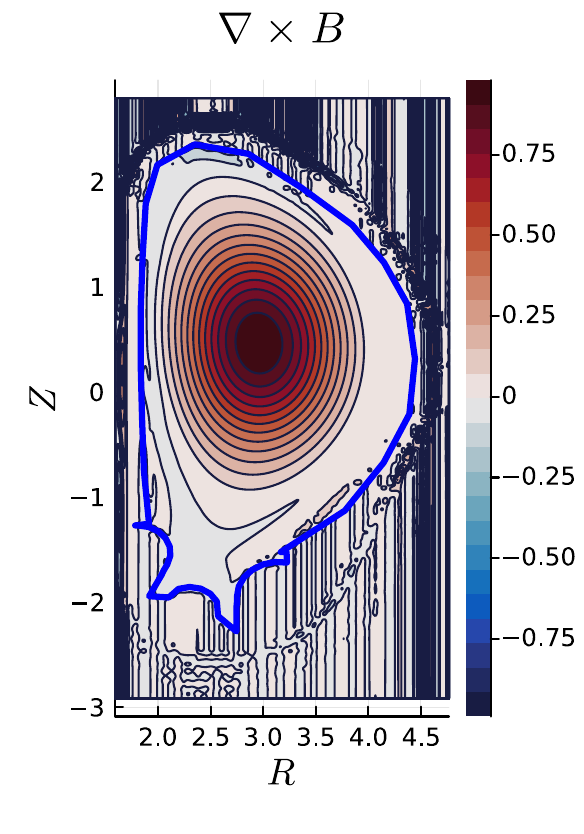}

\caption{\label{fig:fields_mfd} 9th order accurate Hermite interpolation of ITER magnetic fields, reconstructing the poloidal flux function and toroidal current from a mimetic finite difference discretization. The number of grid points in each direction is $N_R = 23$,  $N_Z = 48$. The resulting field approximation is continuous, differentiable and divergence free.
}
\end{figure}

The purpose of this numerical experiment is to show whether higher smoothness improves the simulation of guiding center motion. The fields are computed by a few steps of a mimetic finite difference method solving the quasistatic magnetohydrodynamic equations in an ITER tokamak domain \cite{jorti2023mimetic}, and are then interpolated using Hermite interpolation. The interpolated fields are shown in FIG.~(\ref{fig:fields_mfd}).

As a measure of accuracy, we will look at how well the magnetic moment $\mu$ and toroidal canonical momentum $p_\varphi$ are conserved. In dimensionless variables those become
\begin{equation}
  p_\varphi =\frac 1\omega_c \xi p b_\varphi R - \psi, \ \
  \mu = \frac{(1-\xi^2)p^2 } B.
\end{equation}
% \begin{figure}
%   \centering
% \includegraphics[width=0.7\columnwidth]{img/cons_rk_comparison.pdf}
% \caption{\label{fig:cons_comp} Conservation test. Average, relative consevation error of $p_\varphi$ and $\mu$ as a function of the timestep. Timestep is normalized by the 0.5 ms, being the simulation final time. Solid lines display the conservation errors of Runge-Kutta 4th order method. Dashed lines display the conservation errors of 5th order Runge-Kutta from Dormand-Prince pair. Different colors represent the number of derivatives interpolated with Hermite method for field reconstruction.}
% \end{figure}

We simulate a randomly distributed particle bunch for half a millisecond, which corresponds to \(10^5\) -- \(10^8\) timesteps, and calculate the discrepancy of \(p_\varphi\) and \(\mu\) between the initial time and the final time. The comparison is shown in FIG.~(\ref{fig:cons_comp}). In the plot, we display the relative conservation error of \(p_\varphi\) and \(\mu\) as a function of the timestep averaged over the particle distribution.
 The timestep is normalized by the simulation's final time (0.5 ms). This test has been
performed with different field interpolation variants, varying $m$,
i.e. the number of interpolated derivatives. As we can see, the
5th order Runge-Kutta method cannot achieve the theoretical
convergence rate in the case of an insufficient smoothness ($m < 3$) for the ITER-type finite difference-based magnetic field, i.e. the numerical solution of the quasistatic equations when the timestep is small. In contrast, for the analytic magnetic field described in the previous section, $m=2$ is sufficient. Similar findings hold true for 4th order Runge-Kutta method. The reason for such an inconsistency between numerical and analytic fields is that the interpolation of analytic fields has significantly smaller jumps in high derivatives across the gridlines, although $m=2$ only guaranties $C^2$ continuity of $\psi$, ($C^1$ for the \(B_Z\) component).

Poincaré plots are a practical, effective way of representing complex 3D magnetic fields conveniently. Poincaré plots are generated from the intersection of different magnetic field streamlines with a particular poloidal plane. Provided the streamlines are integrated long enough, it is possible to visualize the shape of the magnetic surfaces. However, the slightest deviation from a divergence-free magnetic field would result in largely perturbed Poincaré plots.

In the case of having to analyze fields arising from finite element approximations that are only weakly divergence free, using the proposed methodology below becomes instrumental. Both to recover strictly solenoidal magnetic fields, and to be able to efficiently and accurately integrate the streamlines required to produce Poincaré plots. Moreover, strongly divergence free approximations also benefit from an increased performance when integrating the streamlines.

In order to reconstruct the fields, we evaluate the solution at each node of a Cartesian grid covering the whole finite element domain. For points outside the original domain, the magnetic field is set to be $\vec{0}$, which would be a bad solution as the field is evaluated close to the boundary, e.g., for diagnostics beyond the separatrix. In the latter case, it would be beneficial to extend the finite element domain beyond the first wall or to extend the Hermite interpolation beyond the finite element boundary by imposing boundary conditions \cite{beznosov2021hermite}. Nevertheless, for flux reconstruction it is important to put a point in the center of the domain as the integration bound (i.e. the point $(R_c, Z_c)$ in \eqref{eq:psi_approx}, \eqref{eq:chi_approx} and perform the line integration in 4 directions rather then 2 directions from the corner of the domain, so that the error from the outside of the domain does not propagate inside the domain through integration. Unfortunately such error propagation will inevitably happen in non-convex regions, i.e. below the X-point, around divertor. Then, we use
\eqref{rec:psi}--\eqref{rec:chi} to find $\psi$ from the
given $\mathbf B$. Once the flux function is known, we are able to efficiently evaluate the magnetic field and produce the Poincaré plots. To this end, we seek the intersection of the integral curves of \(\mathbf{B}\) with a poloidal plane \(\varphi = \varphi_0\). That is, we are looking for the solution curves of
\begin{equation}
  \frac {{\rm d}\mathbf{X}}{{\rm d}t} = \mathbf{B}(\mathbf{X}).
  \label{eq:poincare_plot_orig}
\end{equation}
In a cylindrical coordinate system, taking \(\varphi\) as an independent variable and applying \eqref{eq:int_expr1}, \eqref{eq:int_expr2}, we obtain
\begin{align}
  \frac {{\rm d}R}{{\rm d}\varphi} &= \frac {RB_R}{B_\varphi} \approx \frac{Ru}{v},\\
  \frac {{\rm d}Z}{{\rm d}\varphi} &= \frac {RB_Z}{B_\varphi} \approx \frac{R}{v} \left( -\int_{Z_c}^Z\frac{\partial u}{\partial R}\, {\rm d}Z'
  + w\bigg|_{Z = Z_c}\right) 
  -\frac{1}{v}\int_{Z_c}^Z\frac{\partial v}{\partial \varphi}\, {\rm d}Z',
  \label{eq:poincare_plot}
\end{align}
which we solve for various initial conditions to the end point $\varphi = 2\pi k$ for the $k$-th toroidal transit to generate the Poincaré plots.

In the next subsections, we show two practical examples for which we have used the methodology described above.

\subsubsection{\label{ssec:Linear_piecewise_continuous_fields}Linear piecewise continuous fields}

Let us consider in this experiment a magnetic field resulting from a first order continuous finite element solution. The solution is obtained using the formulation proposed in a previous work \cite{Bonilla_CMAME_2023}. Using this strategy yields magnetic fields that are not strictly solenoidal, since it is only possible to enforce the projection of $\nabla\cdot B$ onto the discrete finite element space to be null.

\begin{figure}[h]
  \begin{subfigure}[t]{0.4\columnwidth}
    \centering
    \includegraphics[height=0.23\textheight]{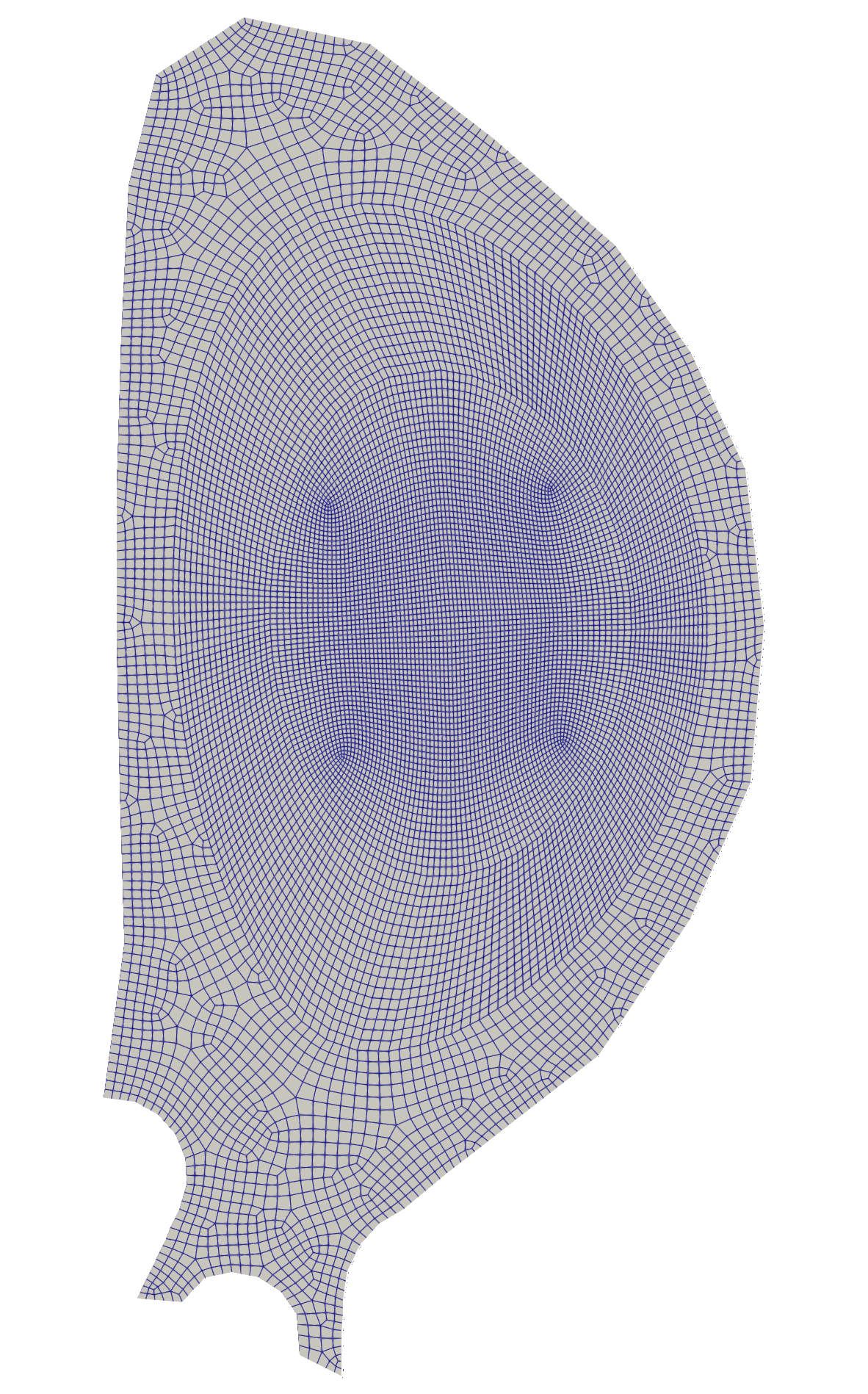}
    \caption{\label{fig:mesh_a} Original finite element mesh where the field is defined.}
  \end{subfigure}
  \begin{subfigure}[t]{0.4\columnwidth}
    \centering
    \includegraphics[height=0.23\textheight]{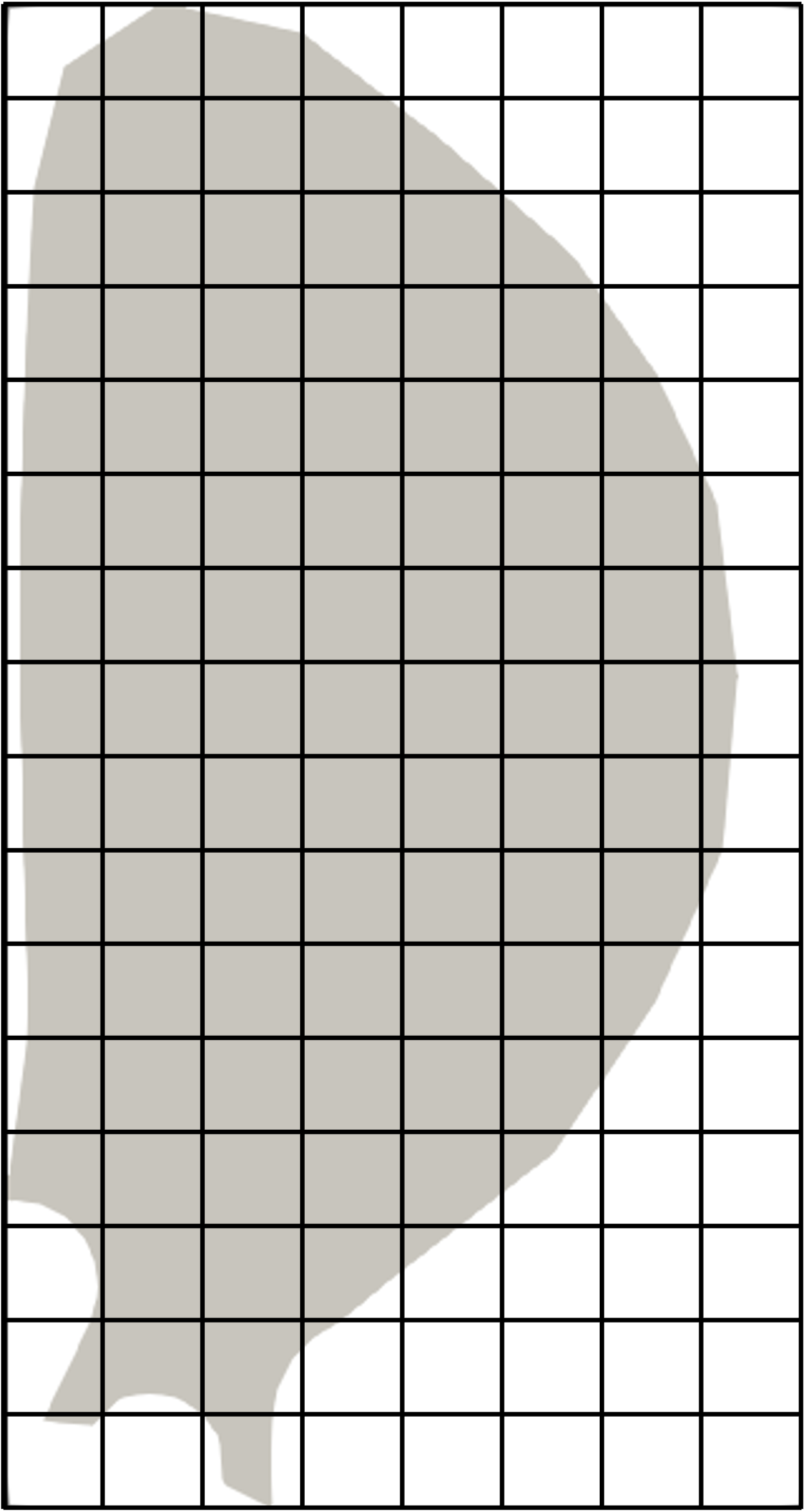}
    \caption{\label{fig:mesh_b} Cartesian mesh where the solution is interpolated for later analysis. Here each box represents a $10\times 10$ mesh block.}
  \end{subfigure}
  \caption{\label{fig:mesh} Meshes used to produce Poincaré plots.}
\end{figure}

We consider a piecewise continuous
linear field given on the mesh depicted in FIG.~\ref{fig:mesh_a}. We
evaluate the solution at each node of an 80 by 160 points Cartesian
grid covering the domain $[1.979,4.3975]\times[-2.205,2.356]$ (see
FIG.~\ref{fig:mesh_b}). Following the previously described procedure,
and integrating equation \eqref{eq:poincare_plot} we are able to efficiently
produce Poincaré plots, which are created using four sets of $21$ seeds and $100$ toroidal transits
for each seed. Each set of seeds is spread equally along lines of length $\{1,1.8,1,2\}$ separated a length 0.2 from the magnetic axis and pitched an angle of $\{0,110,180,252.5\}$ degrees, respectively. Note that all 
lengths are non-dimensionalized with respect to the plasma chamber's
horizontal radius.

\begin{figure}
  \begin{subfigure}[t]{0.19\columnwidth}
    \centering
    \includegraphics[width=\textwidth]{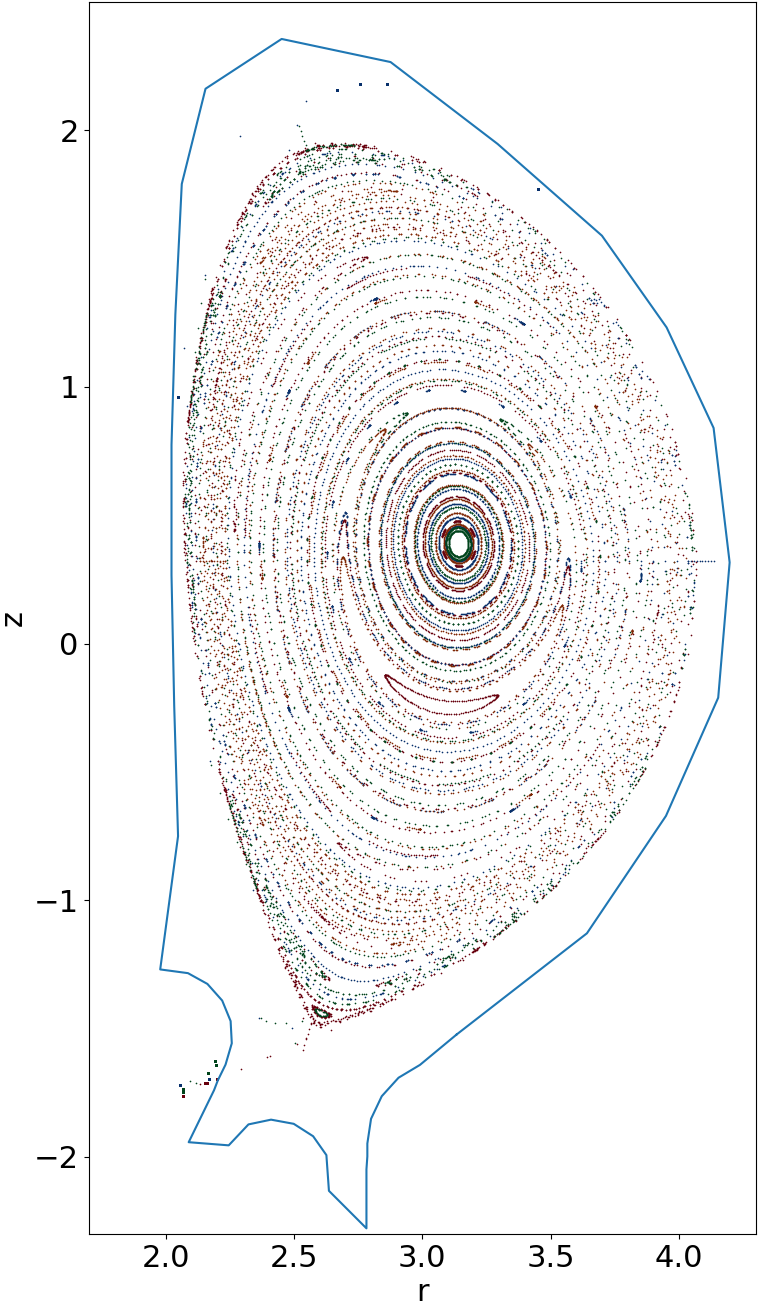}
    \caption{$t=233 \tau_A$.}
  \end{subfigure}
  \begin{subfigure}[t]{0.19\columnwidth}
    \centering
    \includegraphics[width=\textwidth]{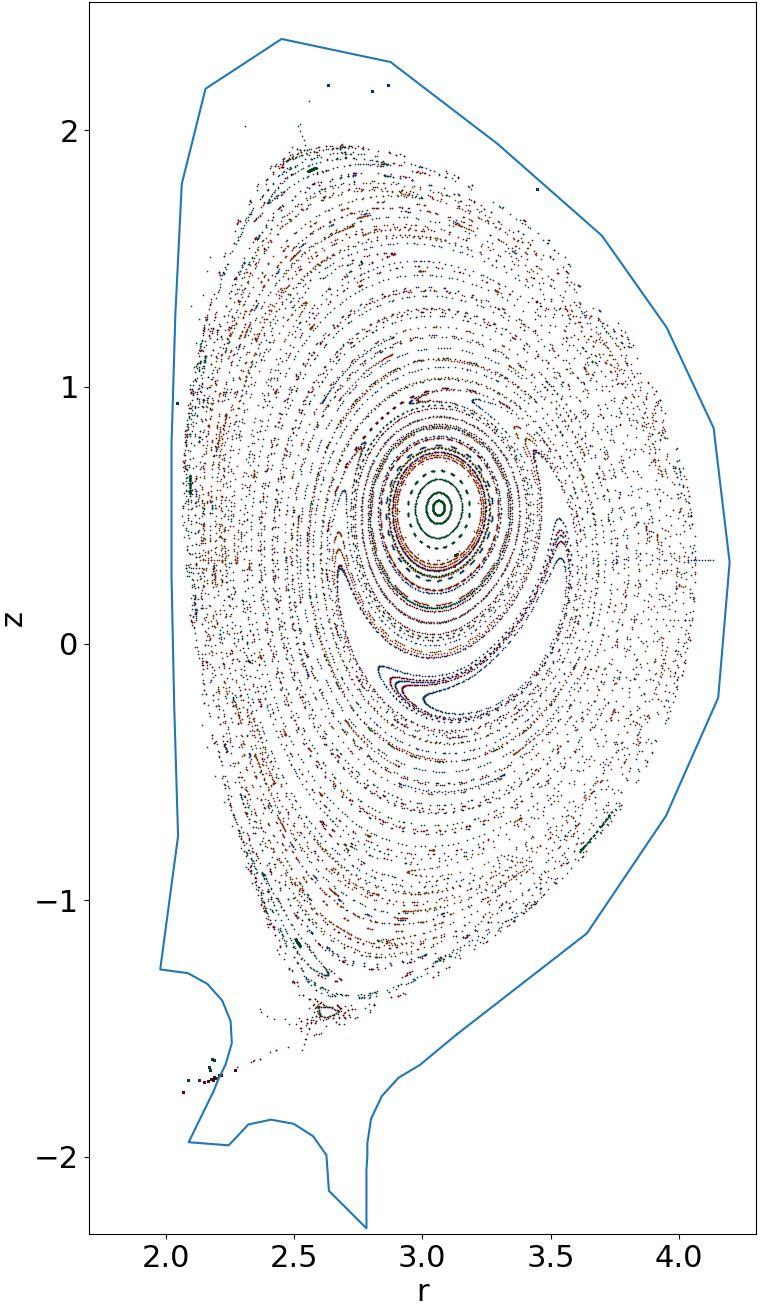}
    \caption{$t=1322 \tau_A$.}
  \end{subfigure}
  \begin{subfigure}[t]{0.19\columnwidth}
    \centering
    \includegraphics[width=\textwidth]{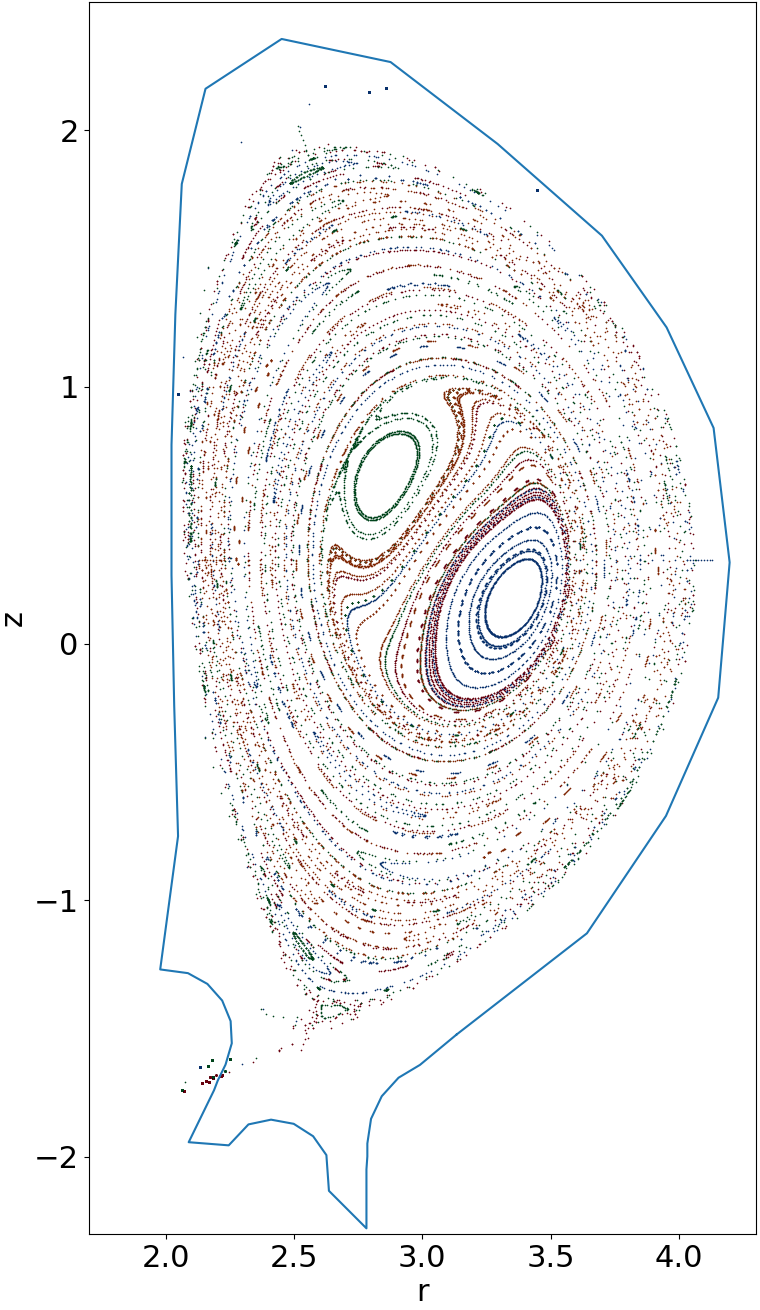}
    \caption{$t=2005 \tau_A$.}
  \end{subfigure}
  \begin{subfigure}[t]{0.19\columnwidth}
    \centering
    \includegraphics[width=\textwidth]{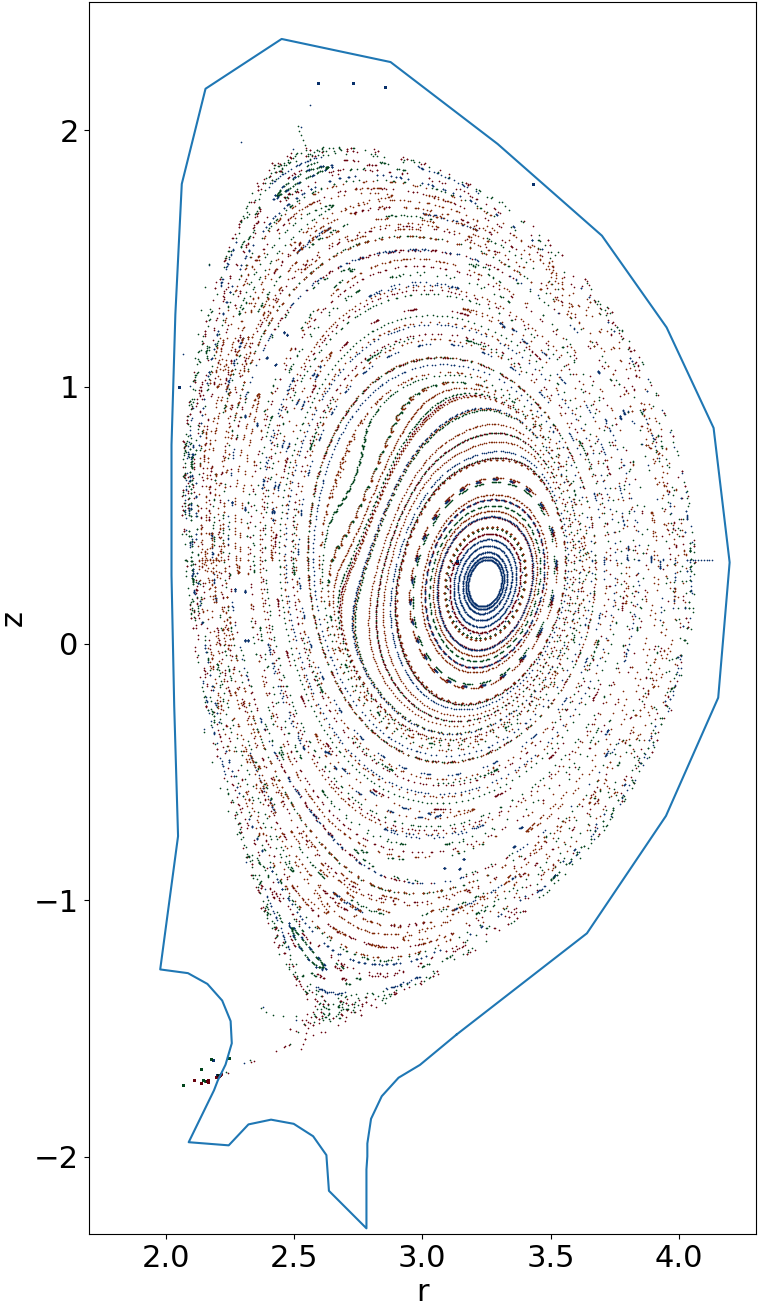}
    \caption{$t=2287 \tau_A$.}
  \end{subfigure}
  \begin{subfigure}[t]{0.19\columnwidth}
    \centering
    \includegraphics[width=\textwidth]{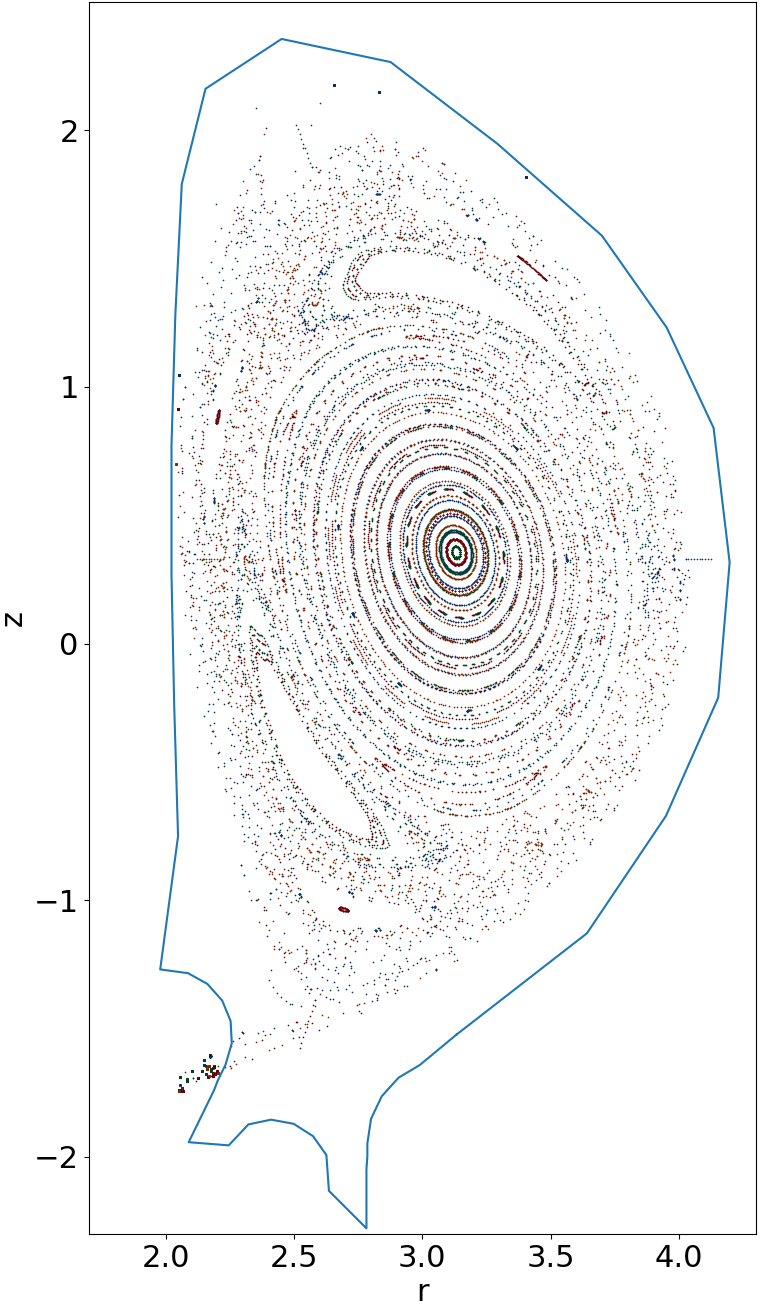}
    \caption{$t=6427 \tau_A$.}
  \end{subfigure}
\caption{\label{fig:poincare} Poincaré plots at different times using the proposed flux reconstruction of an ITER major disruption event~\protect\cite{Bonilla_CMAME_2023}.
}
\end{figure}

FIG.~\ref{fig:poincare} depicts them for different
time steps of a simulation of an ITER major disruption event. Initially, 
a $(1,1)$ kink instability is observed, which excites highly disruptive $(2,1)$ kink mode.

\begin{figure}
  \begin{subfigure}[t]{0.16\columnwidth}
    \centering
    \includegraphics[width=\textwidth]{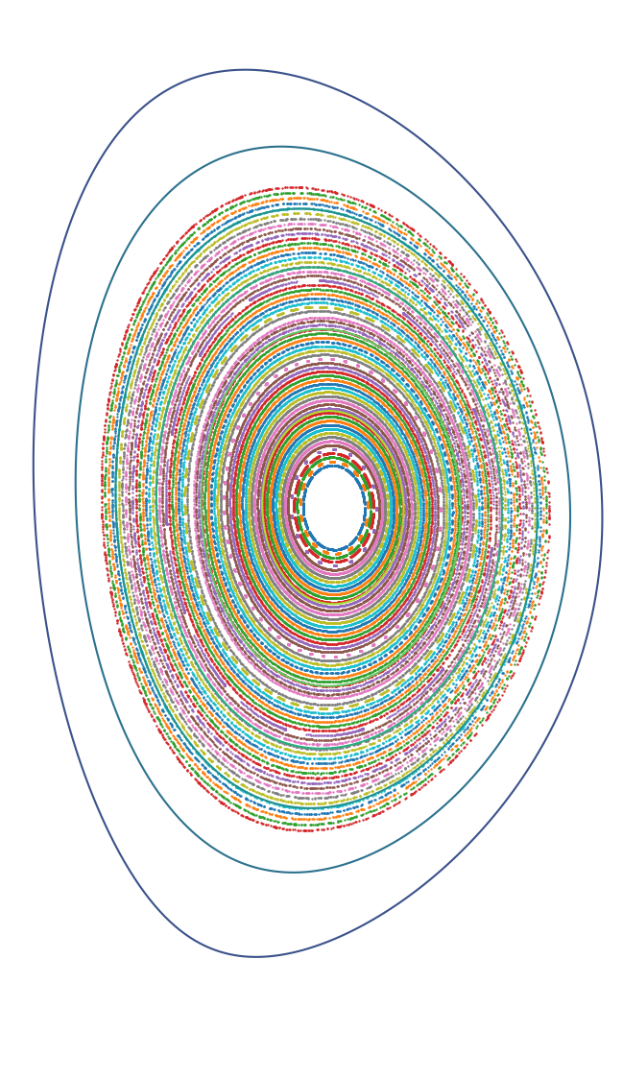}
    \caption{$t=0 \tau_A$.}
  \end{subfigure}
  \hfill
  \begin{subfigure}[t]{0.16\columnwidth}
    \centering
    \includegraphics[width=\textwidth]{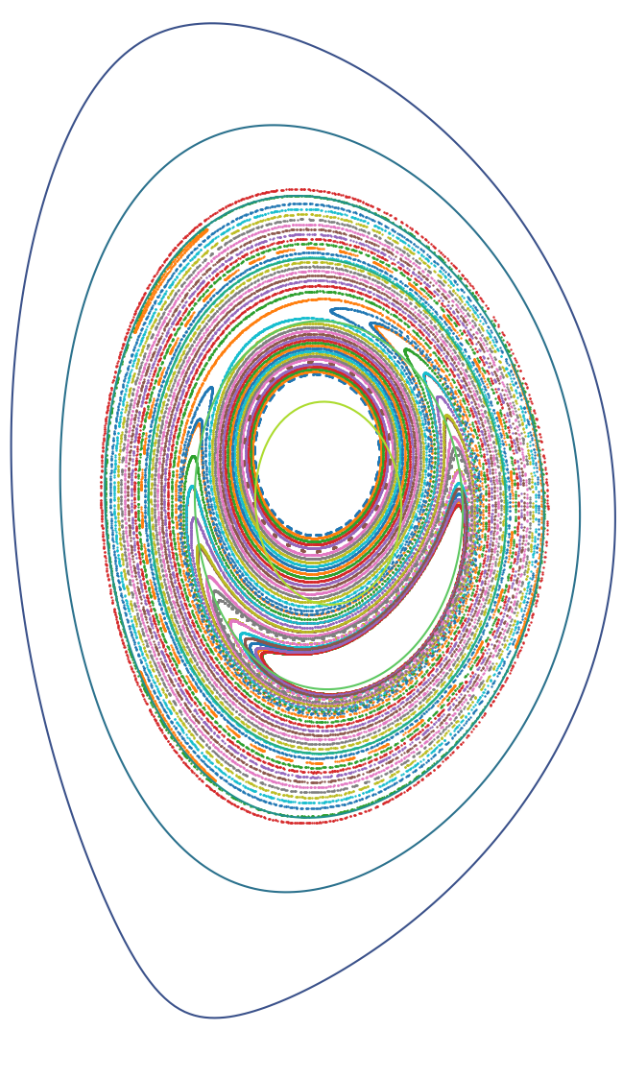}
    \caption{$t=1322 \tau_A$.}
  \end{subfigure}
  \hfill
  \begin{subfigure}[t]{0.16\columnwidth}
    \centering
    \includegraphics[width=\textwidth]{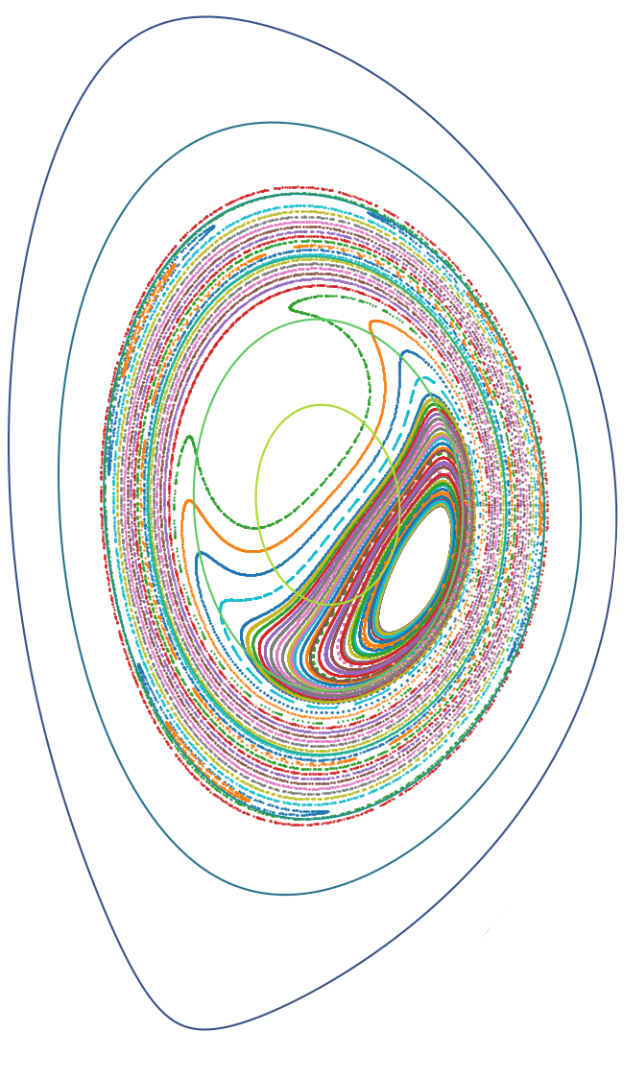}
    \caption{$t=2005 \tau_A$.}
  \end{subfigure}
  \hfill
  \begin{subfigure}[t]{0.16\columnwidth}
    \centering
    \includegraphics[width=\textwidth]{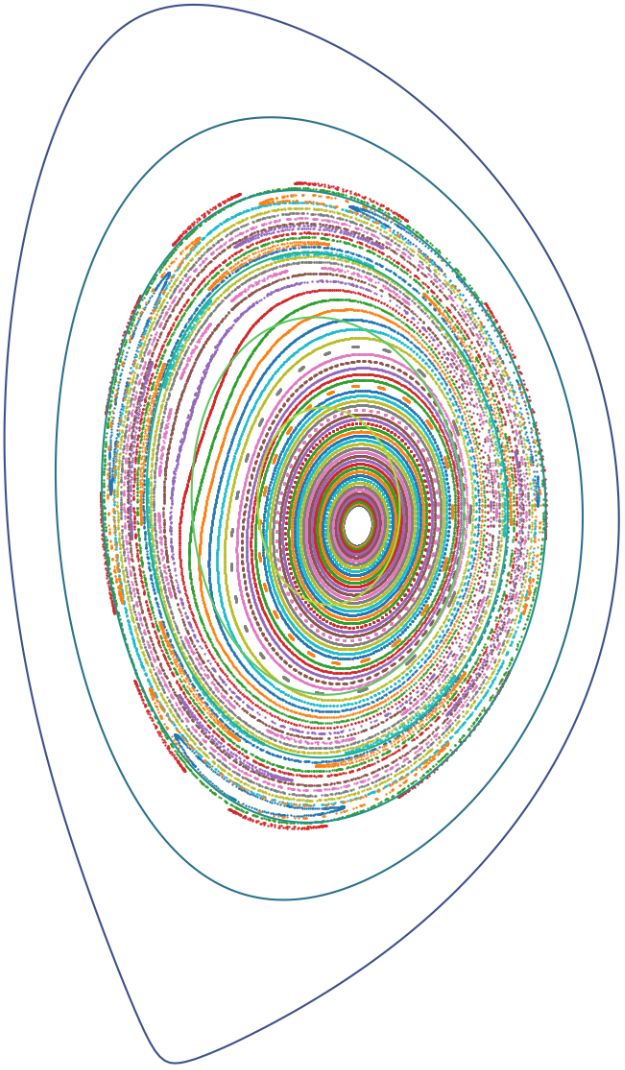}
    \caption{$t=2287 \tau_A$.}
  \end{subfigure}
  \hfill
  \begin{subfigure}[t]{0.16\columnwidth}
    \centering
    \includegraphics[width=\textwidth]{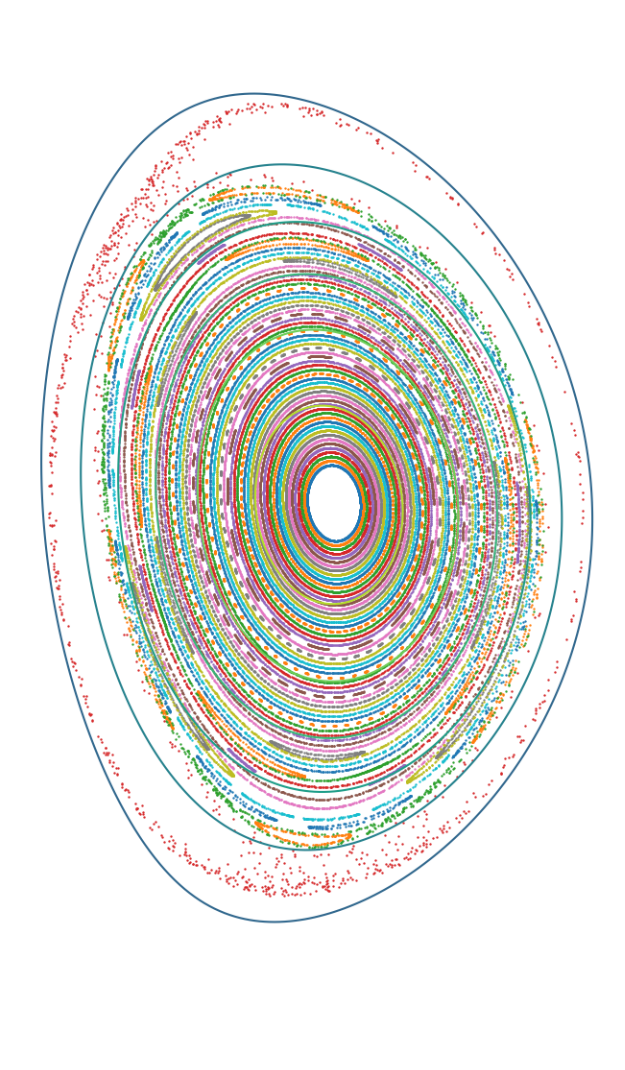}
    \caption{$t=6427 \tau_A$.}
  \end{subfigure}
\caption{\label{fig:psiAndQ} Flux function for the toroidally averaged field and Poincaré plots at different times.}
\end{figure}

Given the toroidally averaged magnetic field is provided, or for axisymmetric fields, it is also relevant to reconstruct and analyze the poloidal flux function. In FIG \ref{fig:psiAndQ}, we show the flux function for the toroidally averaged field and together with the Poincaré plots.

\subsubsection{\label{ssec:Compatible_finite_element_fields}Compatible finite element fields}

\begin{wrapfigure}{r}{8.5cm}
\begin{center}
\vspace{-6mm}
\includegraphics[width=0.48\textwidth]{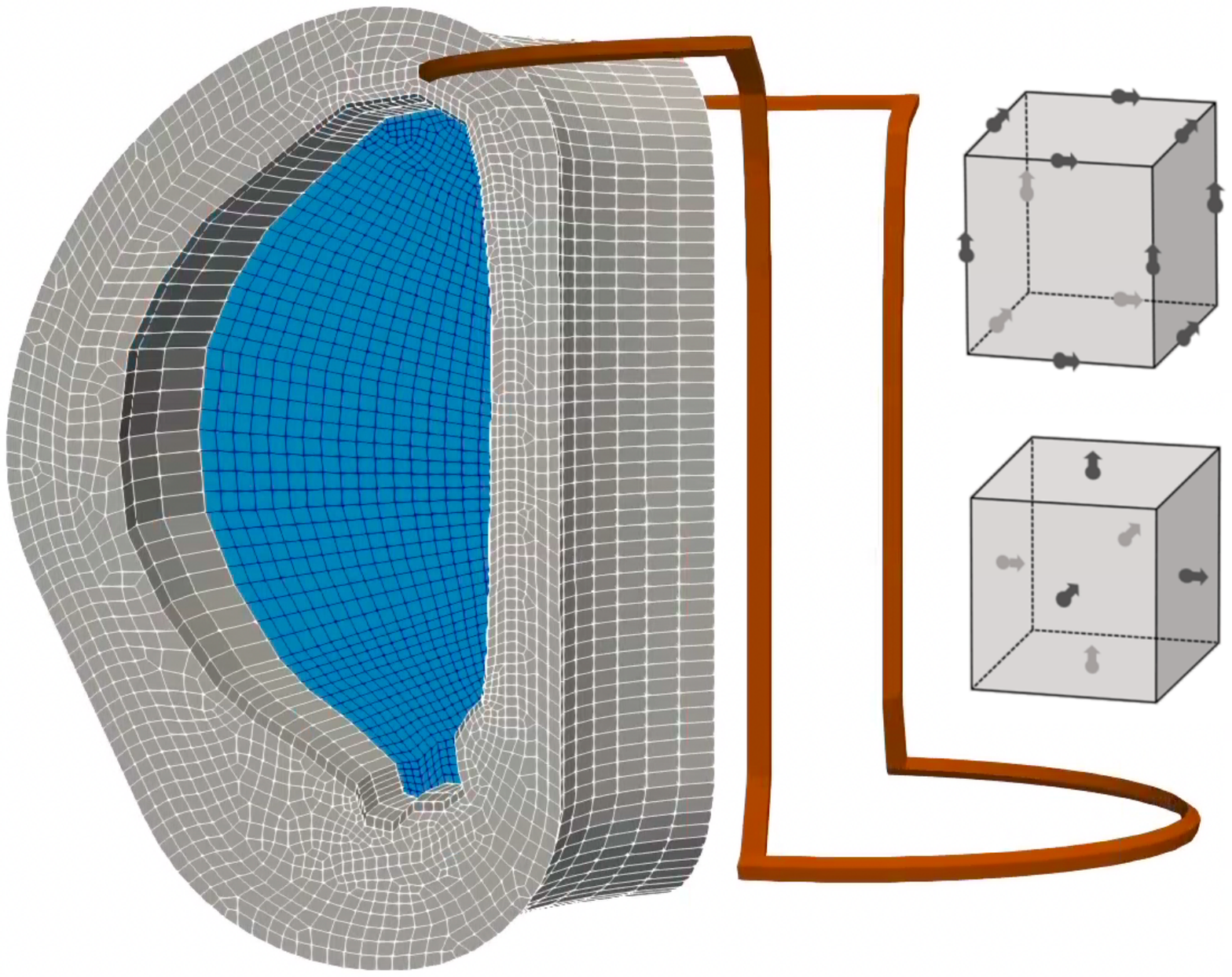}
\end{center}
\vspace{-6mm}
\caption{Mesh and space setup for finite-element based Poincaré
  plots. Left: Mesh cutout including plasma (blue), wall (grey), and
  coil (brown) regions. Right: curl- (top) and div- (bottom)
  conforming hexahedral finite elements for the spaces $Nc_1^e$ and
  $Nc_1^f$, respectively, with degrees of freedom in marked dark
  grey.} \label{MeshSpaceCoil}
\end{wrapfigure}

% Finite element space setup
Next, we test our Poincaré plotting routine on compatible finite elements. 
For this purpose, we first describe our choice of 3D
magnetic field configuration as well as mesh and finite element space
(see Figure \ref{MeshSpaceCoil}).  The magnetic field $\mathbf{B} =
\mathbf{B}_0 + k \mathbf{B}_c$ is given by an MHD-equilibrium field
$\mathbf{B}_0$ -- which corresponds to a typical setup including
nested flux surfaces and one x-point -- together with a
non-axisymmetric perturbation $k \mathbf{B}_c$, for $k \in \{0, 10,
20, ..., 80\}$. The perturbation is given by a magnetic field induced
by a current in a ``helical'' coil embedded inside the reactor
chamber. In practice, such coils may be used to break up closed flux
surfaces in scenarios such as a vertical displacement event (VDE)
after a thermal quench of the plasma.  The 3D stochastic magnetic
fields thus formed can help deplete the high energy runaway electrons.
~\cite{smith_passive_2013} 
%For details on how to construct $\mathbf{B}_c$, see Appendix A.
As in the
mimetic difference case, the magnetic field can be discretized in a
curl- or div-conforming finite element space in order to preserve
(strongly or weakly) the field's zero weak divergence property. Here,
we discretize $\mathbf{B}$ in the lowest order curl-conforming
N\'ed\'elec space of the first kind $Nc_1^e$. Note that this is unlike
the mimetic difference case, which would correspond to
the div-conforming N\'ed\'elec space of the second kind $Nc_1^f$. Finally,
the 3D tokamak mesh $\mathcal{T}_{ext}$ is generated starting from an
irregular 2D poloidal plane mesh $\mathcal{T}_p$, which is then
extruded along the toroidal coordinate $\varphi$ using $40$ poloidal
planes. The two layers including the coil's vertical portion are set
to a thickness of $2\pi/100$ radians; the remaining $38$ layers are
spread evenly over the remaining $2\pi(49/50)$ radians.

% Poincaré plotting setup
In order to compute Poincaré plots for the helically perturbed
magnetic field $\mathbf{B} = \mathbf{B}_0 + k \mathbf{B}_p$, we first
need to process the magnetic field's data. In the case of div-
and curl-conforming finite elements, the data is stored as degrees of
freedom corresponding to integral moments over facets and edges,
respectively. First, we compute the 3D magnetic field's cylindrical
coordinate components $(\mathbf{B}_r, \mathbf{B}_\varphi, \mathbf{B}_z)$
by projecting $\mathbf{B}_i \cdot \mathbf{e}_i$ into a 3D second
polynomial order scalar discontinuous Galerkin space
$dQ_2(\mathcal{T}_{ext})$, for unit vectors $\mathbf{e}_i$ and $i = r,
\varphi, z$. Note that the choice of space ensures that the degrees of
freedom are now associated with coordinates, rather than integral
moments. We then compute the first $l = 5$ Fourier modes in the
toroidal direction for each of these three fields.
This corresponds 
to 
$1 + 2(l-1) = 2l - 1$ mode functions in total, given by the constant
mode together with the sine and cosine field for each of the $(l-1)$
higher modes. In practice, since our mesh is constructed as an
extrusion in the $\varphi$-direction, this can be done using appropriate
summations along $\varphi$ for degrees of freedom sharing the same
$(r,z)$-coordinates. Having integrated over $\varphi$, these 3$\times$(2l
- 1) Fourier fields are then functions in the 2D scalar discontinuous
Galerkin space $dQ_2(\mathcal{T}_p)$. Finally, we interpolate the 2D
Fourier fields into the space $dQ_2(\mathcal{T}_{reg})$, where
$\mathcal{T}_{reg}$ is a regular rectangular 2D mesh corresponding to
a $N_r\times N_z = 40\times80$ grid, and whose underlying rectangular
domain contains the toroidal cross-section domain on which
$\mathcal{T}_p$ is defined. At this point, the fields in
$dQ_2(\mathcal{T}_{reg})$ then have degree of freedom values
associated with the $l$ Fourier modes' values at a regularly spaced
set of coordinates, as required for our Poincaré plotting
routine. Note that overall, this means that the original data for
$\mathbf{B}_0$ is first interpolated from the regular fine MHD
equilibrium code grid onto the semi-irregular coarse finite element
mesh $\mathcal{T}_{ext}$, and then interpolated back to a coarse
regular grid. One may therefore expect a significant interpolation
error that may interfere with the Poincaré plotting routine's
performance.

\begin{figure}[ht]
\begin{center}
\includegraphics[width=.24\textwidth, trim = {10pt 10pt 15pt 15pt}, clip]{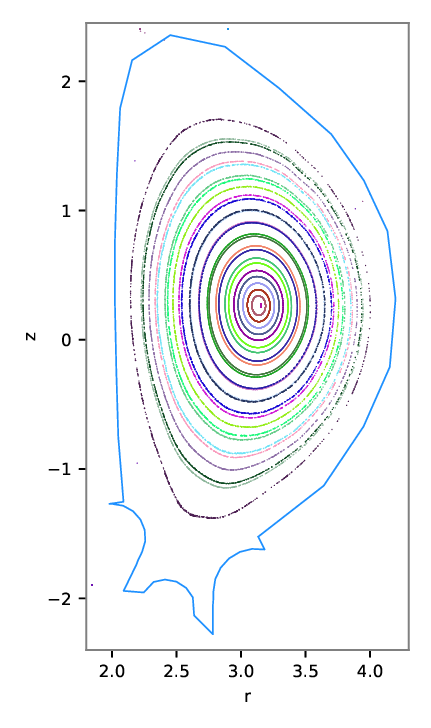}
\includegraphics[width=.24\textwidth, trim = {10pt 10pt 15pt 15pt}, clip]{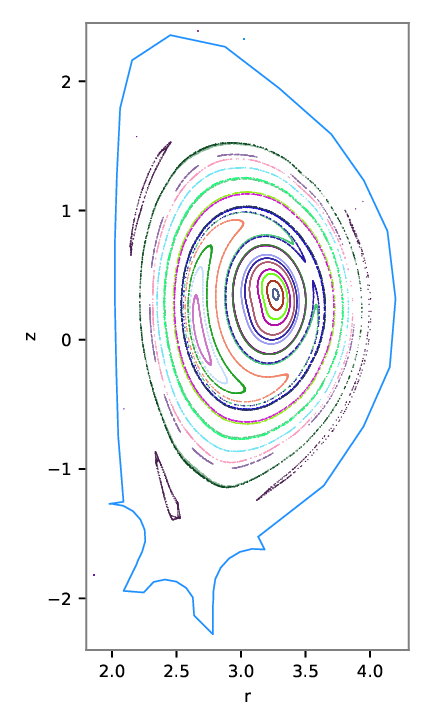}
\includegraphics[width=.24\textwidth, trim = {10pt 10pt 15pt 15pt}, clip]{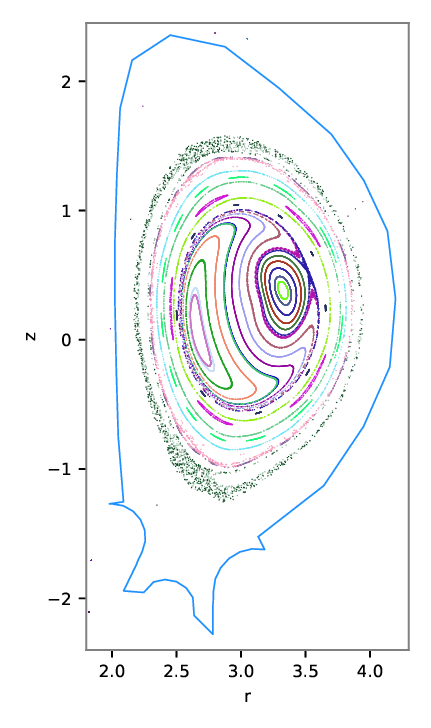}
\includegraphics[width=.24\textwidth, trim = {10pt 10pt 15pt 15pt}, clip]{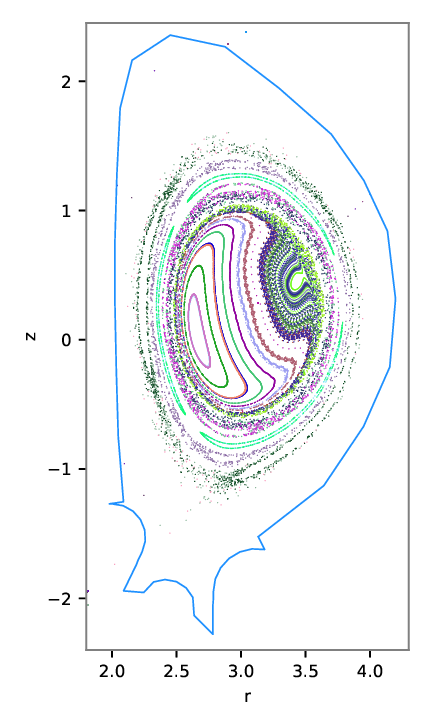}
\includegraphics[width=.24\textwidth, trim = {10pt 10pt 15pt 15pt}, clip]{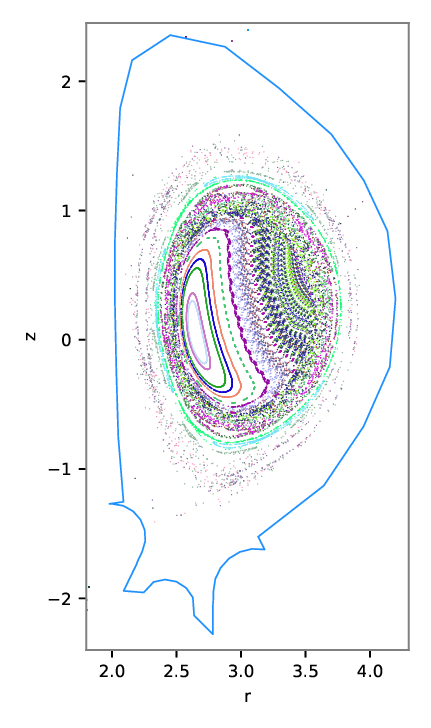}
\includegraphics[width=.24\textwidth, trim = {10pt 10pt 15pt 15pt}, clip]{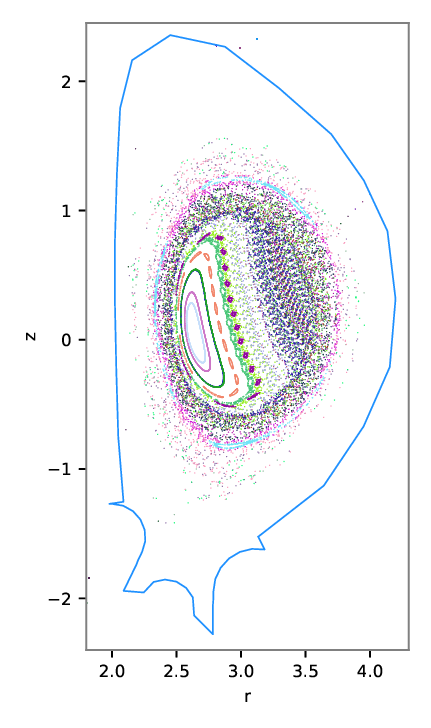}
\includegraphics[width=.24\textwidth, trim = {10pt 10pt 15pt 15pt}, clip]{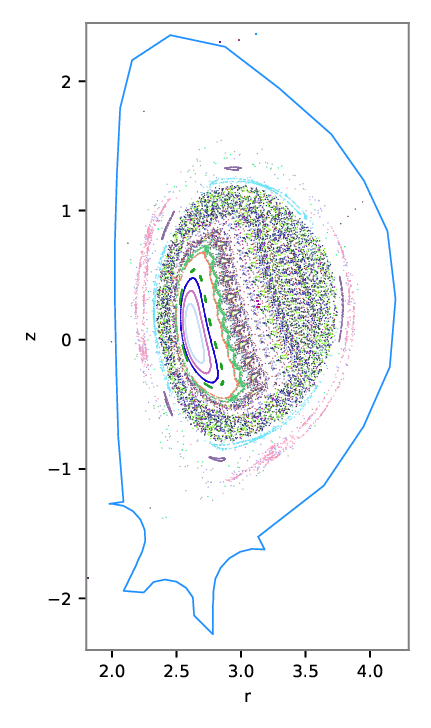}
\includegraphics[width=.24\textwidth, trim = {10pt 10pt 15pt 15pt}, clip]{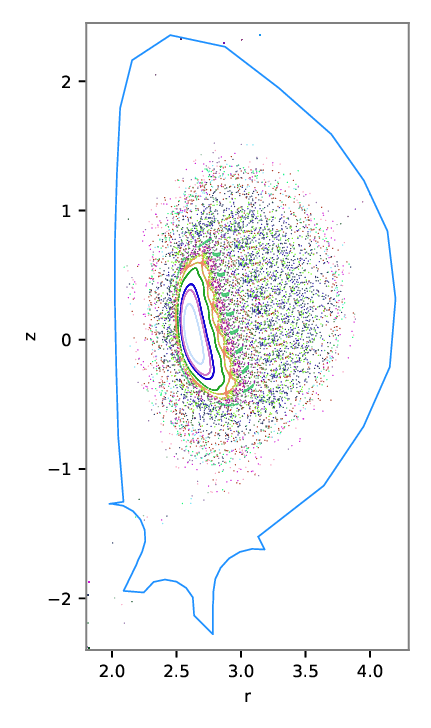}
\end{center}
\vspace{-0mm}
\caption{Poincaré plots for equilibrium magnetic field perturbed by a
  helical coil, for perturbation strengths $k \in \{0, 10, 20, ...,
  80\}$ (left to right, top to bottom). Following a
  non-dimensionalization, the $R$ and $Z$ coordinates are scaled by
  $L_0 = 1/2$.} \label{B_coil_Poincare}
\end{figure}

The Poincaré plots are created using $31$ seeds and $1000$  {toroidal transits}
  for each seed. The seeds
are spread equally along a line of length $2.4$ (noting that the
domain is non-dimensionalized with respect to the plasma chamber's
horizontal radius), containing the magnetic axis of $\mathbf{B}_0 + k
\mathbf{B}_c$, pitched at an angle of $\pi/4$ radians. The resulting
plots are shown in FIG.~\ref{B_coil_Poincare}.

The case $k=0$ corresponds to the equilibrium field $\mathbf{B}_0$ only, and we find that the Poincaré plotting routine recovers the closed flux surface structure very well, in spite of the aforementioned interpolation error. As $k$ increases, we find that as expected, the perturbing magnetic field arising from an asymmetric coil current increasingly breaks up the surfaces, leading to an increasingly stochastic magnetic field profile. Altogether, this demonstrates that our Poincaré plotting routine effectively reproduces the expected magnetic field behavior for lowest order compatible spaces defined on coarse, semi-regular finite element meshes.

%\newpage

\section{\label{sec:Summary}Summary}
We presented a method for reconstructing divergence-free continuous
magnetic fields via vector potentials, with a choice of high order
continuity based on Hermite interpolation. It is designed for
applications when high order adaptive ODE solvers are used to push
particles or integrate the field lines for Poincaré section analysis,
especially if the optimal integration order and a timestep adjustment is
desired. If the input field data is exact, the guaranteed $C(m)$
continuity comes with a $(2m)$--$(2m+2)$ order of accuracy of field
components and poloidal flux function approximation. This was
demonstrated in section~\ref{ssec:Convergence}.

Further, the volume preserving and high order continuity properties of
Hermite interpolation with vector potential reconstruction were tested on the relativistic
Guiding Center Equations. Here, it is crucial to guarantee
conservation of toroidal canonical momentum and magnetic moment, at
least to a small tolerance for accuracy of long term simulations. In
section~\ref{ssec:First_order_mimetic_finite_difference_fields}, it
has been shown that insufficient derivative continuity of the discrete background magnetic field
does degrade the order of Runge-Kutta schemes, especially at low error
tolerances.  We studied this behavior in section~\ref{sec:OrderReduction}, and our findings can be summarized as follows. Consider a
simple ODE $y' = f(y)$ and local truncation error of a 4th order
Runge-Kutta method, in the vicinity of a 3rd derivative discontinuity
\begin{equation}
\varepsilon \sim C\Delta t^5 + \alpha (f^{(3)}_+ - f^{(3)}_-)\Delta t^4,
\end{equation}
where $\alpha$ and $C$ are constants that depend on $f$, and the $+/-$ notation indicates the right and left 3rd derivative of the right-hand side function. In the perfectly smooth case, the second term would cancel out and the local truncation error becomes $O(\Delta t^5)$. Additionally, the method will also achieve a near optimal order of convergence in the discontinuous case, for high error tolerances when the first term dominates and the jump $(f^{(3)}_+ - f^{(3)}_-)$ is small.

Alternatively, Runge-Kutta methods can be used for discontinuous
right-hand sides preserving the optimal order of convergence using
dense output techniques (also based on Hermite interpolation) as
described in \cite{hairer1993}. However, in many particle applications this
may become inefficient, as it requires employing root-finding
mechanisms for each particle individually each time it crosses
an interface. Nevertheless, this would be an interesting method to test in future
work.

\section{\label{sec:Acknowledgement}Acknowledgement}
We thank the U.S. Department of Energy (DOE) Office of Fusion Energy
Sciences and Office of Advanced Scientific Computing Research for
support under the Tokamak Disruption Simulation SciDAC project and the
magnetic fusion theory program, at the Los Alamos National Laboratory (LANL)
under DOE Contract No. 89233218CNA000001.  Additional support to Beznosov
and Tang was provided by the LANL LDRD Program under Explotary
Research project 20250556ER.  This research used resources provided by
the Los Alamos National Laboratory Institutional Computing Program,
which is supported by the U.S. Department of Energy’s National Nuclear
Security Administration under Contract No. 89233218CNA000001, and the
National Energy Research Scientific Computing Center (NERSC), a
U.S. Department of Energy Office of Science User Facility located at
Lawrence Berkeley National Laboratory, operated under Contract
No. DE-AC02-05CH11231 using NERSC award FES-ERCAP0028152.

\bibliography{main}% Produces the bibliography via BibTeX.

\end{document}